%% file: master.tex
\documentclass{article}
\usepackage{times}
\usepackage[hyperindex=true,pageanchor=true,hyperfigures=true,backref=false]{hyperref} 
\usepackage{amsmath}
\usepackage{amssymb}
\usepackage{mathtools}

\usepackage{url}
\usepackage{dsfont} 
\DeclareSymbolFontAlphabet{\Bbb}{AMSb}
\usepackage[T1]{fontenc}

\newlength{\myleftmargin}
\input{macros}

\input{locmacros}

\title{When does a Gaussian process have its paths in a reproducing kernel Hilbert space?}

\author{Ingo Steinwart\footnote{I want to thank Aleksandar Arsenijevic and Daniel Winkle  for helpful discussions. Funded by Deutsche Forschungsgemeinschaft (DFG, German Research Foundation) under Germany’s Excellence Strategy - EXC 2075 – 390740016.}\\
University of Stuttgart\\
Faculty 8: Mathematics and Physics\\
Institute for Stochastics and Applications\\
D-70569 Stuttgart Germany \\
\texttt{\small ingo.steinwart@mathematik.uni-stuttgart.de}
}

\newcommand{\todo}[1]{}

\begin{document}

\maketitle

\begin{abstract}
We investigate for which Gaussian processes there do or do not exist reproducing kernel Hilbert spaces (RKHSs)
that contain almost all of their paths. In particular, we establish a  new result that makes it possible to
 exclude the existence of 
such RKHSs in many cases. Moreover, we  combine this negative result with some known techniques to establish positive results.
Here it turns out that for many classical families of Gaussian processes we can fully characterize for which members of these families 
there exist RKHSs containing the paths. Similar characterizations are obtained for Gaussian processes, for which the RKHSs of their
covariance functions are  Sobolev spaces or Sobolev spaces of mixed smoothness.
\end{abstract}

\textbf{Mathematical Subject Classification (2010).} Primary 60G15, 60G17; Secondary  46E22, 46N30.


\input{intro}

\input{prelims}

\input{results}

\input{sobol}

\input{examples}

\input{proofs}

\bibliographystyle{plain}
\bibliography{../../literatur-db/steinwart-mine,../../literatur-db/steinwart-books,../../literatur-db/steinwart-article,../../literatur-db/steinwart-proc}



%

\end{document}

%% file: macros.tex
\usepackage{amsmath}
\usepackage{amssymb}
\usepackage{amsthm}
\usepackage{color}
\DeclareSymbolFontAlphabet{\Bbb}{AMSb}
\newcommand{\frk}[1]{\mathfrak{#1}}

\newtheorem{theorem}{Theorem}[section]
\newtheorem{lemma}[theorem]{Lemma}
\newtheorem{proposition}[theorem]{Proposition}
\newtheorem{corollary}[theorem]{Corollary}
\newtheorem{definition}[theorem]{Definition}

\newtheorem{example}[theorem]{Example}

\newcommand{\atob}[2]{\emph{#1)} $\Rightarrow$ \emph{#2)}.} 
 
\newcommand{\ada}[1]{\emph{#1).}} 
\newcommand{\adat}[1]{\emph{#1)}} 

\newenvironment{proofof}[1]{\vspace*{1ex}\noindent{\bf Proof of #1:}}{\qed\medskip}

\newlength{\fixboxwidth}
\definecolor{darkgreen}{rgb}{0,0.6,0}

\newcommand{\changed}[1]{\textcolor{black}{#1}}
\newcommand{\changebegin}{\color{black}}
\newcommand{\changeend}{\color{black}}

\newcommand{\scriptonlyraw}[1]{}

\newcommand{\ca}[1]{{\cal #1}}

\newcommand{\eul}{\mathrm{e}}
\newcommand{\imi}{\mathrm{i}}

\newcommand{\mycdot}{\,\cdot\,}

\newcommand{\myqquad}{\qquad\qquad}

\newcommand{\N}{\mathbb{N}}

\newcommand{\R}{\mathbb{R}}
\newcommand{\Rd}{\mathbb{R}^d}

\newcommand{\Cfield}{\mathbb{C}}

\renewcommand{\P}{\mathrm{P}}

\newcommand{\lbd}{\lb^d}

\newcommand{\intd}{\, \mathrm{d}}

\newcommand{\E}{\mathbb{E}}

\renewcommand{\a}{\alpha}
\renewcommand{\b}{\beta}

\renewcommand{\d}{\delta}

\newcommand{\e}{\varepsilon}

\renewcommand{\k}{\kappa}
\newcommand{\lb}{\lambda}

\newcommand{\s}{\sigma}

\newcommand{\om}{\omega}
\newcommand{\Om}{\Omega}

\DeclareMathOperator{\spann}{span}

\DeclareMathOperator{\ran}{ran}
\DeclareMathOperator{\rank}{rank}

\newcommand{\cfp}[1]{\varphi_{#1}}

\DeclareMathOperator{\id}{id}

\newcommand{\snorm}[1]{\Vert #1 \Vert}
\newcommand{\mnorm}[1]{\bigl\Vert \, #1 \, \bigr\Vert}

\newcommand{\inorm}[1]{\Vert #1 \Vert_\infty}

\newcommand{\sborel}{\frk B}

\newcommand{\sLx}[2]{{\ca L_{#1}(#2)}}

\newcommand{\Lx}[2]{{L_{#1}(#2)}}

\newcommand{\skprod}[2]{\langle #1, #2 \rangle}
\newcommand{\skprodb}[2]{\bigl\langle #1, #2 \bigr\rangle}

\newcommand{\tridia}[6]{
\setlength{\unitlength}{1ex}
\begin{picture}(60,20)
\put(20,14){\makebox(9,2)[rt]{$#1$}}
\put(51,14){\makebox(9,2)[lt]{$#2$}}
\put(35.5,2){\makebox(9,2){$#3$}}
\put(30,15){\vector(1,0){20}}
\put(30,14){\vector(1,-1){9}}
\put(41,5){\vector(1,1){9}}
\put(37,16){\makebox(6,2){$#4$}}
\put(27.5,8){\makebox(6,2)[rt]{$#5$}}
\put(47.5,8){\makebox(6,2)[lt]{$#6$}}
\end{picture}}

\newcommand{\quadiass}[8]{
\setlength{\unitlength}{1ex}
\begin{picture}(60,20)
\put(20,15){\makebox(9,2)[rt]{$#1$}}
\put(45,15){\makebox(9,2)[lt]{$#2$}}
\put(20,2){\makebox(9,2)[rt]{$#3$}}
\put(45,2){\makebox(9,2)[lt]{$#4$}}
\put(30,16){\vector(1,0){14}}
\put(28,14){\vector(0,-1){9}}
\put(46,14){\vector(0,-1){9}}
\put(30,3){\vector(1,0){14}}
\put(34,17){\makebox(6,2){$#5$}}
\put(21,8.5){\makebox(6,2)[rt]{$#6$}}
\put(47,8.5){\makebox(6,2)[lt]{$#7$}}
\put(34,0){\makebox(6,2){$#8$}}
\end{picture}}

%% file: locmacros.tex
\newcommand{\sA}{\ca A}

\newcommand{\sB}{\ca B}

\newcommand{\eqclass}[1]{[#1]_\sim}

\newcommand{\sppath}[1]{#1_{\bullet}}

\DeclareMathOperator{\hoel}{H\ddot o l}
\newcommand{\holspace}[2]{\hoel^{#1}(#2)}

\newcommand{\hstar}[1]{H_{#1}^\star}
\newcommand{\kstar}[1]{k_{#1}^\star}
\newcommand{\hdag}[1]{H_{#1}^\dagger}
\newcommand{\htim}[1]{H_{#1}^\times}

\newcommand{\smspac}{\circeq}

%% file: intro.tex
 \section{Introduction}\label{sec:intro}
 
 Given a centered Gaussian process $X:=(X_t)_{t\in T}$ on some (complete) probability space $(\Om,\sA,\P)$,
 a classical task is to describe its path properties by function spaces such as spaces of continuous functions. 
 One particular instance of this general task asks, whether the paths $\sppath X$ of the process are almost surely contained in a suitable   
  reproducing kernel Hilbert space (RKHS) $H$ 
 on $T$. Clearly, the most obvious choice of such a space is the RKHS $H_X$ of the 
  the covariance function $k_X(s,t) := \E(X_s X_{t})$ of $X$. In this case, 
  a classical result, which goes back to Parzen \cite{Parzen63a}, Kallianpur \cite{Kallianpur70a}, and LePage \cite{LePage73a}
 shows that
 \begin{align*}
  \P(\{\sppath X \in H_X\}) = 1
 \end{align*}
if and only if $\dim H_X<\infty$. In other words, in basically no situation of interest is $H_X$ an RKHS containing the paths of $X$.
However, $H_X$ is, of course, by no means the only RKHS on $T$ and beginning with Driscoll \cite{Driscoll73a}, and 
 subsequently Fortet \cite{Fortet73a},
more general $H$ have been considered. 
  The most  general investigation  in this direction has been conducted by 
  Luki\'c and Beder \cite{LuBe01a} and since their results provide a basis of our analysis, we 
  quickly recall their findings. To this end, 
%
%
%
%
%
%
we need the following notion.

\begin{definition}\label{def:nuc-dom}
 Let $H_1$ and $H_2$ be RKHSs on $T$. Then we write $H_1 \ll H_2$ and speak of nuclear dominance of $H_2$ over $H_1$, if 
  $H_1\subset H_2$ and 
 the resulting embedding map $J:H_1\to H_2$ is a Hilbert-Schmidt operator. 
\end{definition}

Note that 
formally, \cite{LuBe01a} defines nuclear dominance differently, see the discussion following  Theorem \ref{thm:dominance} below, but Lemma \ref{lem:comput-dominance} shows that 
our definition is actually equivalent to the one used in \cite{LuBe01a}. With these preparations, \cite[Theorem 7.4]{LuBe01a}
provides the following characterization for a \emph{fixed} RKHS $H$.

\begin{theorem}\label{thm:LuBe01a}
 Let $X:=(X_t)_{t\in T}$  be a centered Gaussian process on a complete probability space
 $(\Om,\sA,\P)$ and $H$ be an RKHS on $T$. Then   the following statements hold true:
 \begin{enumerate}
  \item If $H_X\ll H$, then there exists a version $Y$ of $X$ with $\P(\{\sppath Y \in H\}) = 1$.
 \item If $H_X\not \ll H$, then for all versions $Y$ of $X$  we have $\P(\{\sppath Y \in H\}) = 0$.
 \end{enumerate}
\end{theorem}

Here we note that strictly speaking  \cite[Theorem 7.4]{LuBe01a} only shows \adat {ii} for $Y=X$. However, if we have an 
arbitrary version $Y$ of $X$, then $Y$ is a centered Gaussian process 
with $k_Y = k_X$, which in turn gives   $H_Y=H_X$. Applying the original result 
\cite[Theorem 7.4]{LuBe01a}  to $Y$ then yields  \adat {ii}.

Informally  speaking,  given a 
 fixed ``candidate'' RKHS $H$, Theorem \ref{thm:LuBe01a} shows that we have  either  
 \begin{align*}
 \P(\{\sppath X \in H\}) = 1 \myqquad \mbox{ or } \myqquad \P(\{\sppath X \in H\}) = 0\, ,
 \end{align*}
 and in addition both situations 
 are characterized by $H_X\ll H$, respectively $H_X\not \ll H$.
 In this sense, it provides a ``test'' for each RKHS $H$ on $T$.
However, it does not tell us how to \emph{find} an $H$ with $H_X\ll H$, nor whether  $H_X\ll H$   is even possible for 
suitable RKHSs $H$.

In this respect we note that
positive, constructive
results can sometimes be derived if we already have some knowledge on the paths of $X$. 
To illustrate this, 
 let $E$ be a  
Banach space   of functions (BSF) on $T$, that is, $E$ is a Banach space consisting   
of functions $f:T\to \R$ such that the point evaluations $\d_t:E\to \R$ given by $\d_t(f) := f(t)$ for all $t\in T$ and $f\in E$, are continuous. 
If $(\Om,\sA,\P)$ is complete and we   know that there exists a version $Y$ of $X$ with $\P(\{\sppath Y \in E\}) = 1$,
then for every ``surrounding'' RKHS
$H$ on $T$, that is $E\subset H$, we obviously  also have
$\P(\{\sppath Y \in H\}) = 1$.
Unfortunately, however, this technique is limited to ``sufficiently benign'' $E$, since
there are various important BSFs $E$, for which 
no such surrounding  $H$ exists. For example, for  
the somewhat natural candidate $E=C(T)$, where $C(T)$ denotes 
the space of   continuous functions $T\to \R$ on a compact metric space $T$, 
we cannot find a surrounding $H$ as soon as $T$ is not countable, see e.g.~\cite{Steinwart24a}.
Moreover,  \cite{ScSt25a} has recently shown that  even for many spaces $E$ of smoother functions no such surrounding $H$ with bounded kernel may exist.
For example,  
 if $\holspace \b {[0,1]}$  denotes  the space of all  $\b$-H\"older continuous functions $f:[0,1]\to \R$, then there exists an RKHS $H$ with bounded kernel and $\holspace \b {[0,1]}\subset H$
if $\b>1/2$, while there is no such RKHS if $\b<1/2$, see \cite[Section 4.5,  Theorem 22, and Lemma 24, respectively Theorem 12]{ScSt25a}.
In particular, the approach of taking a surrounding RKHS does not work for the Wiener process, while for the 
fractional Brownian motions with Hurst index $> 1/2$ 
it does work, see also Example \ref{ex:fbm}.
In summary, the approach of using a surrounding RKHS can only be used if we already know that the paths are contained in a sufficiently small BSF $E$.

A different approach
for constructing  an RKHS $H$  with $\P(\{\sppath Y \in H\}) = 1$ for some suitable version $Y$ of $X$
has been taken in \cite{Steinwart19a}. In a nutshell, this approach is based 
on a refined analysis of the Karhunen-Lo\`eve expansion for $X$. Unfortunately, this  
approach, as well as its recent refinement in \cite{Karvonen23a}, again fails to produce  positive results in many situations. In particular,   it does not 
provide an RKHS containing the paths of the Wiener process.

The example of the Wiener process suggests that  for some Gaussian processes $X$ it may be \emph{impossible} 
to find an RKHS containing the paths of $X$. 
To the best of our knowledge, however, such 
general 
impossibility results, even for concrete processes,  have not been established   yet, despite the fact that 
Gaussian processes with paths in an RKHS have  attracted interests in diverse areas of applied 
mathematics, see e.g.~\cite{PiWuLiMuWo07a,FlSeCuFi16a,Gualtierotti15,TrGi24a} and the references mentioned in \cite{Karvonen23a}.
In particular,  Theorem \ref{thm:LuBe01a} cannot be directly used, since we would have to check  $H_X\ll H$ for \emph{all}
RKHS $H$ on $T$, or at least 
for a sufficiently rich class of such RKHSs. Similarly, a failure of the described two approaches for finding an $H$ containing 
the paths of the considered process does, of course,  not imply a general impossibility result. 
%
 
 The main goal of this paper is close this gap by  
 \emph{characterizing} when 
 for a Gaussian process $X$ there exists an RKHS $H$ and a version $Y$ of $X$ with  $\P(\{\sppath Y \in H\}) = 1$.
 To this end, we will first present 
  a  \emph{necessary} condition
 for the existence of $H$ with $H_X\ll H$. Combining this necessary condition with the 
 two approaches discussed above for finding $Y$ and $H$ with 
 $\P(\{\sppath Y \in H\}) = 1$ we can then characterize for various families of Gaussian processes, when such $Y$ and 
 $H$ exist.
%
%


The rest of this  paper is organized as follows: In Section \ref{sec:prelims} we collect the notions and notations that are necessary 
for understanding our results. The general theory is then presented in 
Section \ref{sec:results},
and in  Section  \ref{sec:sobol-rkhs}  we apply our methods to Gaussian processes $X:=(X_t)_{t\in T}$ 
whose RKHSs $H_X$ are essentially equal to a fractional Sobolev space on $T\subset \R^d$.
In addition, the case of Sobolev spaces with dominating mixed smoothness is also investigated.
Moreover, some concrete  examples are presented in 
Section \ref{sec:examples}. Finally, all proofs as well as some additional material needed for these proofs
can be found in Section \ref{sec:proofs}.

%
%
%
%

%% file: prelims.tex
\section{Preliminaries}\label{sec:prelims}

For sequences $(a_n)_{n\geq 1}$ and $(b_n)_{n\geq 1} \subset (0,\infty)$ we write $a_n\sim b_n$ if there exist
constants $c_1,c_2>0$ with $c_1 a_n \leq b_n\leq c_2 a_n$ for all $n\geq 1$. 

If $E$ and $F$ are Banach spaces,
 we write $E\hookrightarrow F$, if $E\subset F$ and the resulting embedding map is continuous. 
In addition, we write $E\smspac F$ if $E\subset F\subset E$, and we write  
 $E\cong F$ if $E\smspac F$ and the norms of $E$ and $F$ are equivalent. In other words, we have 
 $E\cong F$ if and only if $E\hookrightarrow F\hookrightarrow E$.
Obviously $E\cong F$ implies $E\smspac F$, and if $E$ and $F$ are BSFs, then the converse implication is also true 
by a simple application of the closed graph theorem, see e.g.~\cite[Lemma 23]{ScSt25a}. For RKHSs $H_1$ and $H_2$ on $T$
we thus have $H_1 \smspac H_2$ if and only if $H_1\cong H_2$.

\changed{Given a non-empty set $T$, we write $\ell_\infty(T)$ for the space of all bounded functions $T\to \R$ equipped with the  supremums norm. Clearly, $\ell_\infty(T)$ is a BSF. Moreover,  if  $(T,\sB)$ is a  measurable space,}
 we write $\sLx 0\sB$ for the space of measurable functions $T\to \R$, where we equip 
$\R$ with the Borel $\s$-algebra $\sborel$. Analogously, $\sLx \infty\sB$ stands for the Banach space of all bounded, measurable functions $T\to \R$
 equipped with the supremums-norm $\inorm\cdot$. 
Furthermore, if $\nu$ is some measure on $(T,\sB)$, then $\Lx 2\nu$ denotes the usual Hilbert space of $\nu$-equivalence classes $\eqclass f$ of square $\nu$-integrable functions $f$,
and $\Lx \infty \nu$ stands for the Banach space of $\nu$-equivalence classes $\eqclass f$ of measurable, $\nu$-almost surely bounded functions $f$.
Note that the map
\begin{align}\label{eq:Inu}
I_{\nu}: \sLx \infty \sB  &\to \Lx \infty \nu \\ \nonumber
f &\mapsto \eqclass f \, ,
\end{align}
is linear and bounded with 
$\inorm{\eqclass f} \leq \inorm f$ for all $f\in \sLx \infty \sB$, where $\inorm{\eqclass f}$ denotes
the usual essential supremum norm of $\eqclass f$ in $\Lx \infty \nu$.
Finally, if $\lb^d$ is the Lebesgue measure on $\Rd$ and 
$T\subset \Rd$ is measurable with $\lb^d(T)>0$, then we often write $\Lx 2 T:= \Lx 2 {\lb^d_{|T}}$, where $\lb^d_{|T}$ is the Lebesgue measure on $T$.
Analogously, we write $\sLx \infty T := \sLx \infty{\sborel^d_{|T}}$, where $\sborel^d_{|T}$ denotes the restriction of the  Borel $\s$-algebra 
$\sborel^d$ on $\Rd$ to $T$. 

Let us now assume that we have an RKHS $H$ on $T$ with kernel $k$, where we refer to e.g.~\cite{Aronszajn50a,BeTA04}, 
and \cite[Chapter 4]{StCh08}
for some background information on these spaces.
Then 
$k$ is bounded if and only if all functions in   $H$ are bounded, see e.g.~\cite[Lemma 4.23]{StCh08}.
Moreover, if $\sB$ is some $\s$-algebra on $T$ for which $k$ is
$\sB\otimes \sB$-measurable, then $H\subset \sLx 0 \sB$, see e.g.~\cite[Lemma 4.24]{StCh08}, and 
the converse implication is true for separable $H$, see e.g.~\cite[Lemmas 4.24 and 4.25]{StCh08}.
Finally, let us assume that we have a finite measure space $(T,\sB,\nu)$ and an RKHS $H$ on $T$ with $H\subset \sLx \infty \sB$.
As mentioned above this yields $H\hookrightarrow \sLx \infty \sB$, 
and therefore the linear map 
\begin{align}\label{def:Ikn-prelim}
I_{k,\nu} :H&\to \Lx 2 \nu\\ \nonumber
h&\mapsto \eqclass h
\end{align}
is well-defined and continuous. 
Finally,  we say that $\nu$ is $H$-positive, if $I_{k,\nu}$ is injective, that is, $\nu(\{h\neq 0\})>0$ for all $h\in H$ with $h\neq 0$.
For example, if $H$ consists of continuous functions and $\nu(U)>0$ for all non-empty, open $U\subset T$, then $\nu$ is $H$-positive.
On the other hand, if e.g.~$\nu$ is an empirical measure based upon $n$ points
and $\dim H=\infty$, then $I_{k,\nu}$ cannot be injective.

Recall that 
a probability measure $\P$ on $\R$ is a Gaussian measure, if 
there are $\mu\in \R$ and $\s\in [0,\infty)$ such that 
the characteristic function of $\P$ is given by 
\begin{align*}
\cfp \P(t) =\eul^{\imi  t\mu} \cdot \exp\biggl( -\frac{ \s^2 t^2}{2} \biggr) 
\, , \qquad \qquad t\in \R.
\end{align*}
Clearly, if $\s>0$, then $\P$ is a normal distribution with expectation $\mu$
and variance $\s^2$.
In the case $\s = 0$, the  probability measure $\P$ is a point mass on $\mu$, that is $\P(\{\mu\}) = 1$.

Now let $X= (X_t)_{t\in T}$ be an $\R$-valued, stochastic process  on some probability space $(\Om,\sA,\P)$. As usual, 
we call the map $\sppath X(\om): t\mapsto X_t(\om)$, where $\om\in \Om$ is fixed,  a path of $X$. Moreover, $X$ 
is called a Gaussian process, if for all
$t_1,\dots,t_n\in T$ and $\a_1,\dots,\a_n\in \R$, the distribution of $\a_1 X_{t_1}+\dots+\a_n X_{t_n}$ is a Gaussian measure on $\R$.
Furthermore, the Gaussian process is centered, if $\E X_t = 0$ for all $t\in T$. In this case, its covariance function $k_X:T\times T\to \R$ is given by 
\begin{align*}
k_X (t_1,t_2) := \E (X_{t_1} X_{t_2})\, , \myqquad t_1,t_2\in T.
\end{align*}
It is well known and easy to verify that $k_X$ is a symmetric and positive definite function, and thus a (reproducing) kernel. 
We write $H_X$ for the corresponding   RKHS of $k_X$ and we further use the 
short hand   $I_{X,\nu} := I_{k_X,\nu}$.


%% file: results.tex
\section{General Results}\label{sec:results}

In this section we present and investigate necessary and sufficient conditions 
for the existence of RKHSs containing the paths of a given Gaussian process. 

\changebegin

To this end, we first need a notion that controls the growth behavior of the functions
in a given RKHS.

\begin{definition}
Let $H$ be an RKHS on $T$. Then we say that a function $b:T\to [0,\infty)$ is an envelope for $H$ if, for all $h\in H$ and $t\in T$, we have 
\begin{align}\label{eq:envelope}
|h(t)| \leq \snorm h_H \cdot b(t)\,.
\end{align}
\end{definition}

At first glance,  Condition \eqref{eq:envelope} my look restrictive. However, given an arbitrary RKHS $H$ on $T$, a simple estimate yields 
\begin{align}\label{eq:all-envelope}
|h(t)| = \bigl|  \langle h, k(\mycdot, t)  \rangle_H  \bigr| 
\leq \snorm h_H \cdot \snorm {k(\mycdot, t)}_H
= \snorm h_H \cdot \sqrt{k(t,t)} 
\end{align}
for all $h\in H$ and $t\in T$. In other words, every RKHS has an envelope function, namely the diagonal function
$t\mapsto  \sqrt{k(t,t)}$. In general, however, the 
measurability of this diagonal function does not follow from the measurability of all $h\in H$, and this missing implication may 
cause some technical issues if we would directly work with the diagonal function.
Moreover, the diagonal function
is closely related to the considered RKHS, whereas choosing an envelope function is to some extend independent of a specific RKHS. For example, 
by considering all RKHSs $H$ on $T$ having the constant envelope function $b  = 1$, we actually work with all RKHSs $H$
satisfying $\snorm{\id : H\to \ell_\infty(T)} \leq 1$. For more information 
on the relation between  diagonal and envelope functions we refer to Lemma \ref{lem:envelope-vs-kernel}.

Let us now assume that we have a centered Gaussian process 
$X:= (X_t)_{t\in T}$ whose RKHS satisfies  $H_X\subset \sLx 0 \sB$, where $\sB$ is some $\s$-algebra on $T$.
If $b:T\to [0,\infty)$ is a measurable envelope function for $H_X$
then Lemma \ref{lem:Iknu-new} shows that 
$I_{X,\nu}:H_X\to \Lx 2 \nu$ defined in \eqref{def:Ikn-prelim} 
is a Hilbert-Schmidt operator for all  measures $\nu$ on $(T,\sB)$ with $b\in \sLx 2 \nu$.
%
Our first theorem shows that $I_{X,\nu}$ is even nuclear if the paths of 
a suitable version $Y$ of $X$ are contained in some RKHS $H$ with envelope function $b$.

\begin{theorem}\label{thm:general-result-envelope}
Let $(\Om,\sA,\P)$ be a complete probability space and $X:= (X_t)_{t\in T}$ be a centered Gaussian process over  $(\Om,\sA,\P)$
with  RKHS $H_X$. Moreover, let 
$H$ be an RKHS on $T$ 
and $Y$ be a version  of $X$ with 
\begin{align*}
\P(\{\sppath Y \in H\}) >0\, .
\end{align*}
Then, for   all  $\s$-algebras $\sB$ on $T$ with $H_X\subset \sLx 0 \sB$, all measurable envelope functions $b:T\to [0,\infty)$ of $H$,
%
%
and all 
 measures $\nu$ on $(T,\sB)$ with $b\in \sLx 2 \nu$, the operator $I_{X,\nu}:H_X\to \Lx 2 \nu$ is well-defined and nuclear.
\end{theorem}

Note that there is a trade-off between the considered class of RKHS $H$ in Theorem \ref{thm:general-result-envelope} and 
the strength of the derived necessary condition. Indeed, if we have two measurable envelope functions $b_1$ and $b_2$ with $b_1\leq b_2$,
then every RKHS $H$ having envelope $b_1$ also has envelope $b_2$. In contrast, not every measure $\nu$ with $b_1\in \sLx 2 \nu$
also needs to satisfy $b_2\in \sLx 2 \nu$. In other words, by increasing the envelope function we potentially shrink the class of considered 
operators $I_{X,\nu}:H_X\to \Lx 2 \nu$. An extreme case in this trade-off are the constant envelope functions
 that describe RKHSs containing bounded functions. The following
corollary presents the results of Theorem \ref{thm:general-result-envelope}  for such envelopes.

\begin{corollary}\label{cor:general-result} \changeend
Let $(\Om,\sA,\P)$ be a complete probability space, $X:= (X_t)_{t\in T}$ be a centered Gaussian process over  $(\Om,\sA,\P)$
with  RKHS $H_X$, and $\sB$ be a  $\s$-algebra on $T$ with \changebegin 
$H_X\subset \sLx 0 \sB$\changeend. 
If there exist an RKHS $H$ on $T$ with $H\subset \changebegin \ell_\infty(T) \changeend$
and a version $Y$ of $X$ such that 
\begin{align*}
\P(\{\sppath Y \in H\}) >0
\end{align*}
then, for   all 
finite measures $\nu$ on $(T,\sB)$, the operator $I_{X,\nu}:H_X\to \Lx 2 \nu$ is 
\changebegin well-defined \changeend and  nuclear.
\end{corollary}
 
\changeend

Note that  there is a conceptual difference between 
Theorems \ref{thm:LuBe01a} \changebegin and \ref{thm:general-result-envelope}, respectively Corollary \ref{cor:general-result}. \changeend 
Indeed, Theorem \ref{thm:LuBe01a} provides a ``test'' for every single candidate RKHS $H$ as discussed in the introduction. 
In contrast, Theorem \ref{thm:general-result-envelope} and \changebegin its  Corollary \ref{cor:general-result} \changeend  provide a general impossibility result: \changebegin For example, if \changeend 
for a given  $H_X\subset \sLx \infty\sB$ we find
a finite measure $\nu$ on $(T,\sB)$ for which $I_{X,\nu}:H_X\to \Lx 2 \nu$ is \emph{not} nuclear, then 
there exist  no RKHS $H\subset \sLx \infty \sB$
and no version $Y$ of $X$ with 
\begin{align*}
\P(\{\sppath Y \in H\}) >0 \, .
\end{align*}
While at first glance finding such a $\nu$ may seem to be a hard task, we will see in Sections   \ref{sec:sobol-rkhs} and 
\ref{sec:examples} that for e.g.~bounded, measurable $T\subset \R^d$ \changebegin with non-empty interior \changeend
the Lebesgue measure is typically a sufficient choice. In fact, we will leverage there from well-established 
results on the ``embedding''  $I_{X,\nu}:H_X\to \Lx 2 T$ if $H_X$ belongs to particular classes of RKHSs such as fractional 
Sobolev spaces. 

\changebegin

We finish this discussion with another corollary showing 
that a certain ``localization'' reduces the problem to 
bounded RKHS $H$  even if we are interested in more general RKHSs.

\begin{corollary}\label{cor:restrict-to-bounded}
Let $(\Om,\sA,\P)$ be a complete probability space and $X:= (X_t)_{t\in T}$ be a centered Gaussian process over  $(\Om,\sA,\P)$
with  RKHS $H_X$.
Moreover, let 
$H$ be an RKHS on $T$ 
and $Y$ be a version  of $X$ with 
%
\begin{align*}
\P(\{\sppath Y \in H\}) >0 \, .
\end{align*}
Then, for   all  $\s$-algebras $\sB$ on $T$ with $H_X\subset \sLx 0 \sB$, all measurable envelope functions $b:T\to [0,\infty)$ of $H$, all $M>0$ with $T^* := \{t\in T: b(t) \leq M\}\neq \emptyset$, and 
all 
finite measures $\nu$ on $(T^*, \sB_{|T^*})$, the operator 
\begin{align*}
I_{X_{|T^*}, \nu}: H_{X_{|T^*}} \to \Lx 2 \nu
\end{align*}
is nuclear, where 
  $ H_{X_{|T^*}}\subset \sLx \infty {\sB_{|T^*}}$ denotes the RKHS of the restricted process 
$X_{|T^*} := (X_t)_{t\in T^*}$.
\end{corollary}

Note that on $T^*$, every function $f\in H$ is bounded, namely we have
$|f(t)| \leq \snorm f_H \cdot M$ for all $t\in T^*$. In other words, the restriction of $H$ onto $T^*$
consists of bounded, measurable functions, and therefore we are actually in the situation of 
Corollary \ref{cor:general-result}. The proof of Corollary \ref{cor:restrict-to-bounded}
uses this observation.

\changeend

In some other situations, however, the RKHS $H_X$ and the operator $I_{X,\nu}:H_X\to \Lx 2 \nu$ may not be
as amenable as desired. 
In the following, we will thus present two equivalent characterizations of $I_{X,\nu}:H_X\to \Lx 2 \nu$ being nuclear, which 
sometimes make it easier 
to work with this property.

We will begin with a simple yet useful characterization that makes it possible to   alter the space $H_X$ by some finite dimensional space.
To this end, we need the following definition.

\begin{definition}
Let $H_1$ and $H_2$ be RKHSs on $T$. Then we say that $H_1$ and $H_2$ are essentially equivalent,  if there exist
RKHSs $\hstar 1$ and $\hstar 2$ on $T$ with $\dim \hstar 1 < \infty$, $\dim \hstar 2<\infty$, and 
\begin{align}\label{eq:ess-equiv}
H_1 + \hstar 1 \smspac H_2 + \hstar 2\, .
\end{align}
\end{definition}

It is straightforward to check that the notion of being essentially equivalent
gives an equivalence relation on the set of all RKHSs on $T$, but in the following we do not need this.
In addition, another simple calculation shows that 
if $H_1$ and $H_2$ are essentially equivalent, then both $H_1$ and $H_2$ are also
essentially equivalent to $H_1+H_2$.
Furthermore, RKHSs $H_1$ and $H_2$ with $H_1\smspac H_2$ are obviously also essentially equivalent.
Finally,
the following lemma shows that in some situations the essential equivalence can be easily verified with 
the help of the kernels of the involved RKHSs.

\begin{lemma}\label{lem:kernel-diff}
Let $H_1$ and $H_2$ be RKHSs on $T$ with kernels $k_1$ and $k_2$, and  
$\kstar 1$ and $\kstar 2$ be kernels on $T$ with finite dimensional RKHSs. 
Then $H_1$ and $H_2$ are essentially equivalent, if we have 
\begin{align*}
k_1-k_2 =  \kstar 2 - \kstar 1 \, .
\end{align*}
\end{lemma}

Now, the following result shows that $I_{X,\nu}:H_X\to \Lx 2 \nu$ is  nuclear if and only if $I_{H,\nu}:H\to \Lx 2 \nu$
is nuclear for all RKHS $H$ that are essentially equivalent to $H_X$. This will make it possible to 
switch from $H_X$ to an essentially equivalent 
$H$ for which  $I_{H,\nu}:H\to \Lx 2 \nu$ is better understood.

\begin{theorem}\label{thm:ess-eq-rkhs}
Let  $(T,\sB,\nu)$ be a measure space and 
$H_1$ and $H_2$ be essentially equivalent  RKHSs on $T$ with \changebegin
$H_1\subset \sLx 0 \sB$ and $H_1\subset \sLx 0 \sB$.
If, in addition, $b\in \sLx 2 \nu$ is an envelope function for both 
$H_1$ and $H_2$, \changeend
then $I_{H_1,\nu}:H_1\to \Lx 2 \nu$ is nuclear if and only if $I_{H_2,\nu}:H_2\to \Lx 2 \nu$
is nuclear.
\end{theorem}

Theorem \ref{thm:ess-eq-rkhs} can be used to establish non-existence results with the help of 
Theorem \ref{thm:general-result-envelope} \changebegin and its Corollary \ref{cor:general-result}\changeend, see e.g.~Example \ref{ex:bb}.
Combining the
  next result with Theorem \ref{thm:LuBe01a} shows that in some situations
the  essential equivalence of RKHSs can also be used to find RKHSs containing the paths of a given Gaussian process, 
see e.g.~Theorem \ref{thm:sobol-main}.

\begin{theorem}\label{thm:constr-essent-equiv-dom-rkhs}
 Let 
 $X:=(X_t)_{t\in T}$  be a centered Gaussian process on a probability space  $(\Om,\sA,\P)$ 
 and  $H$ be an RKHS
 that is essentially equivalent to $H_X$. If there exists an RKHS $H^\dagger$ on $T$ with $H\ll H^\dagger$, then there also
 exists an RKHS $H_X^\dagger$ on $T$ that is essentially equivalent to $H^\dagger$ and satisfies both $H_X^\dagger \subset H_X + H^\dagger$
 and 
 \begin{align*}
 H_X \ll H_X^\dagger\, .
 \end{align*}
\end{theorem}

 Our next goal is to present another equivalent way to 
express that $I_{X,\nu}:H_X\to \Lx 2 \nu$ is nuclear as well as to \emph{construct} an $H$ with $H_X\ll H$ in some situations.
To this end, \changebegin  we assume that $b\in \sLx 2 \nu$ is an envelope function for $H_X \subset \sLx 0 \sB$,
where we note that 
in the situation of Theorem \ref{thm:general-result-envelope},
we do have $H_X\subset H$
by Theorem \ref{thm:LuBe01a}, and therefore the envelope function of $H$ in Theorem 
\ref{thm:general-result-envelope} is, modulo a constant, also an envelope function for $H_X$.
In this sense, the assumption of $b$ being an envelope for $H_X$ fits to 
our previous results.   \changeend
Let us now
 consider the operator
$T_{X,\nu} := I_{X,\nu}\circ I_{X,\nu}^*:\Lx 2 \nu\to \Lx 2 \nu$, where $I_{X,\nu}^*:\Lx 2 \nu\to H_X$ denotes the adjoint of $I_{X,\nu}$. Then
some elementary calculations, see Lemma \ref{lem:Iknu-new}, show that $T_{X,\nu}$ is given by  
\begin{align*}
T_{X,\nu} f =  \int_T k_X(\mycdot, t) f(t) \intd \nu(t)\, , \myqquad f\in \Lx 2 \nu,
\end{align*}
where the right-hand side is viewed as a $\nu$-equivalence class in $\Lx 2 \nu$. By construction, $T_{X,\nu}$ is positive and self-adjoint,
and it is also compact, since $I_{X,\nu}$ is  Hilbert-Schmidt by our assumption
\changebegin  of $b\in \sLx 2 \nu$ being an envelope for $H_X$, \changeend
see again Lemma \ref{lem:Iknu-new}.
Consequently, the spectral theorem, see Section \ref{subsec:hs+nuc-op}, shows that
$T_{X,\nu}$ has  at most countably many non-zero eigenvalues $(\mu_i(T_{X,\nu}))_{i\in I}$, where we have set either $I:=\{1,\dots,n\}$ or $I:=\N$, 
and where we have included the geometric multiplicities. Moreover, these eigenvalues are positive and 
there exists an orthonormal system (ONS) $(\eqclass{f_i})_{i\in I}\subset \Lx 2 \nu$ such that each $\eqclass{f_i}$
is an eigenfunction 
of $T_{X,\nu}$ with eigenvalue $\mu_i := \mu_i(T_{X,\nu})$.
We now define 
 \begin{align}\label{eq:def-ei}
 e_i := \mu_i^{-1} I_{X,\nu}^* \eqclass{f_i} \, , \myqquad i\in I.
 \end{align}
This gives $e_i\in H_X$ and $\eqclass{e_i} = I_{X,\nu} e_i = \mu_i^{-1} T_{X,\nu} \eqclass{f_i} = \eqclass{f_i}$ for all $i\in I$.
In other words, we may choose representatives $e_i$ of $\eqclass{f_i}$ that are contained 
in $H$. Finally,  $T_{X,\nu}$ is nuclear since 
$I_{X,\nu}$ is  Hilbert-Schmidt, see Lemma \ref{lem:hs-vs-nuc}, and therefore $(\mu_i(T_{X,\nu}))_{i\in I}$ is summable.

Now recall from \cite{StSc12a} that in some cases the eigenfunctions and -values can be used to construct RKHSs surrounding $H_X$. 
To briefly describe this construction, let us fix a $\b>0$ for which we have 
\begin{align}\label{eq:beta-power-kx}
\sum_{i\in I}  \mu_i^\b \, e_i^2(t) < \infty\, , \myqquad t\in T.
\end{align}
Then \cite[Lemma 2.6]{StSc12a} shows that 
\begin{align}\label{eq:beta-power-kx-full}
k_X^{\b} (t,t') := \sum_{i\in I} \mu_i^\b e_i(t) e_i(t')\, , \myqquad t,t'\in T 
\end{align}
defines a kernel on $T$ whose RKHS is given by 
\begin{align*}
H_X^{\b} := \biggl\{ \sum_{i\in I} a_i  \mu_i^{\b/2}   e_i : (a_i)_{i\in I}\in \ell_2(I)   \biggr\}  \, , 
\end{align*}
where the series converges pointwise and $\ell_2(I)$ denotes the space of all 2-summable families 
$(a_i)_{i\in I} \subset \R$. Moreover, $(\mu_i^{\b/2}   e_i)_{i\in I}$ is an ONB of $H_X^{\b}$ and we have $H_X^\b\subset \sLx 0 \sB$.
Obviously, Condition \eqref{eq:beta-power-kx} is monotone in the sense that if it holds for some $\b$ then it also holds for all exponents 
$\a>\b$. 
In addition, for such pairs we have $H_X^\a \subset H_X^\b$.
Finally, if $\nu$ is $H$-positive, then $(\sqrt{\mu_i} e_i)_{i\in I}$ is an ONB of $H_X$, see 
\cite[Theorem 3.1]{StSc12a}, and thus we have $H_X^1 = H_X$. Consequently, if in this case we have found a $\b\in (0,1)$ with
\eqref{eq:beta-power-kx} then 
$H_X \subset H_X^\b$.

With these preparations we can now present our last result of this section, where we note that \emph{ii)} essentially simplifies 
similar results of \cite[Section 5]{Steinwart19a}.

\begin{theorem}\label{thm:ew-char}
Let $(\Om,\sA,\P)$ be a probability space, \changebegin
$(T,\sB, \nu)$ be a  measure space and
$X:= (X_t)_{t\in T}$ be a centered Gaussian process over $(\Om,\sA,\P)$ such that
$b\in \sLx 2 \nu$ is an envelope  for $H_X$ and $H_X\subset \sLx 0\sB$. \changeend
Moreover, let $(\mu_i)_{i\in I}$ be the non-zero eigenvalues of 
$T_{X,\nu}$ and $(e_i)_{i\in I}\subset H_X$ be a family such that $(\eqclass{e_i})_{i\in I}\subset \Lx 2 \nu$ is an ONS in $\Lx 2 \nu$
and each $\eqclass{e_i}$
is an eigenfunction 
of $T_{X,\nu}$ to the eigenvalue $\mu_i$.
Then the following statements hold true:
\begin{enumerate}
\item  The operator $I_{X,\nu}:H_X\to \Lx 2 \nu$ is nuclear if and only if 
\begin{align}\label{eq:double-nuc-int-op}
\sum_{i\in I} \sqrt{\mu_i} < \infty.
\end{align}
\item Assume that  $\nu$ is $H_X$-positive and we have a  $\b\in (0,1)$ satisfying \eqref{eq:beta-power-kx}. Then 
$H_X\ll H_X^{\b}$ if and only if 
\begin{align}\label{eq:nuc-int-op-dominance}
\sum_{i\in I}  \mu_i^{1-\b} < \infty.
\end{align}
In this case, there exists a version $Y$ of $X$  with $\P(\{\sppath Y \in H_X^\b\}) =1$.
\end{enumerate}
\end{theorem}

To illustrate Theorem \ref{thm:ew-char} we assume for a moment that   $\nu$ is $H_X$-positive and $\sup_{i\in I}\inorm{e_i} < \infty$.
Then, if  \eqref{eq:double-nuc-int-op} is \emph{not} satisfied,   Theorem \ref{thm:ew-char} shows that 
 $I_{X,\nu}:H_X\to \Lx 2 \nu$ is not nuclear, and therefore \changebegin Corollary \ref{cor:general-result} \changeend shows that there exist no RKHS $H\subset \sLx \infty \sB$
and no version $Y$ of $X$ with 
\begin{align*}
\P(\{\sppath Y \in H\}) >0 \, .
\end{align*}
Conversely, if \eqref{eq:double-nuc-int-op} is satisfied, then  both \eqref{eq:beta-power-kx} and \eqref{eq:nuc-int-op-dominance} hold true 
for $\b := 1/2$, and hence
there exists a version $Y$ of $X$ with $\P(\{\sppath Y \in H_X^{1/2}\}) =1$.
Furthermore, note that in this case $k_X^{1/2}$ is bounded, and hence we have $H_X^{1/2} \subset \sLx \infty \sB$.
In other words, 
if we have   an $H_X$-positive, finite measure $\nu$ with $\sup_{i\in I}\inorm{e_i} < \infty$, then 
 \eqref{eq:double-nuc-int-op}   \emph{characterizes} the existence of an RKHS $H\subset \sLx \infty \sB$
containing the paths of a suitable version of $X$. We refer to Example \ref{ex:abel} for such a situation.

Finally, recall from \cite[Theorem 4.6]{StSc12a} that the spaces $H_X^\b$ can in some sense 
be identified with the 
interpolation spaces $[\Lx 2 \nu, [H_X]_\sim]_{\b,2}$ of the real 
method.  In particular, if under some mild technical assumptions we have 
$[\Lx 2 \nu, [H_X]_\sim]_{1/2,2}\hookrightarrow \Lx \infty \nu$, then $H_X^{1/2}$ exists and is an RKHS version of  
$[\Lx 2 \nu, [H_X]_\sim]_{1/2,2}$, see \cite[Theorem 5.5]{Steinwart19a}, 
\changebegin or its
  recent refinement \cite[Theorem 3.4]{BiStXXa}, \changeend
in combination with 
\cite[(36) and Theorem 4.6]{StSc12a},
 while  \cite[Theorem 5.3]{StSc12a} ensures \eqref{eq:double-nuc-int-op}.
In such situations it may thus be possible to construct $H_X^{1/2}$ modulo equivalent norms \emph{without} 
using the eigenfunctions and -values.
 In this respect note that  for certain kernels $k_X$ \cite{FlSeCuFi16a} presents another way to find $H_X^{1/2}$ without 
 the eigenfunctions and -values. 
 
 Finally, \cite[Section 1]{Steinwart19a} lists various articles in which the eigenfunctions and -values 
 are computed for specific Gaussian processes. For such processes, we can then apply Theorem \ref{thm:ew-char} 
 directly.

%% file: sobol.tex
\section{Results for RKHSs of Sobolev Type}\label{sec:sobol-rkhs}

In this section we investigate how our main results can be used to investigate Gaussian processes whose RKHSs are essentially equivalent 
to a fractional Sobolev space, or equal a Sobolev space of dominating mixed smoothness. 
To this end,  recall that for  $s\in (0,\infty)$, the fractional Sobolev spaces, or Bessel potential spaces, are defined by
\begin{align*}
 H^s(\R^d) := \bigl\{ f\in \Lx 2 \Rd: (1+\snorm\cdot^2)^{s/2}  \hat f \in \Lx 2 \Rd   \bigr\}\, ,
\end{align*}
where
$\hat f$ is the Fourier transformation of $f$, see e.g.~\cite[Section 3.1]{DNPaVa12a}  or \cite[Section 2.2.2]{Triebel83}.
By setting 
\begin{align*}
 \snorm f_{H^s(\R^d) } := \mnorm{(1+\snorm\cdot^2)^{s/2}  \hat f }_{\Lx 2 \Rd} 
\end{align*}
the space $H^s(\R^d)$ becomes a Hilbert space and for $m\in \N$ we have $H^m(\Rd)  \cong W^m_2(\Rd)$, where 
$W^m_2(\Rd)$ denotes the classical Sobolev space with smoothness $m$, see e.g.~\cite[(3) in Chapter 2.5.6]{Triebel83} or \cite[Chapter 7, Paragraph 7.62]{AdFo03}. Moreover, the definition of the fractional Sobolev spaces immediately ensures 
 $H^s(\R^d) \subset  H^r(\R^d)$ for all $0<r<s <\infty$, and  for $s>d/2$ we have 
\begin{align}\label{eq:bessel-rkhs-short} 
 H^s(\R^d)  \hookrightarrow  C_b(\Rd)
\end{align}
by Sobolev's embedding theorem, see the beginning of Section \ref{sec:proofs-hoel+sobol} for details and references.
In other words, in the case $s>d/2$ the space $H^s(\R^d)$ can be viewed as an RKHS with bounded and continuous kernel.
Finally, for bounded, measurable $T\subset \Rd$ and $s>d/2$ we define 
$H^s(T) := \{f_{|T} :    f\in H^s(\Rd)\}$ and equip this space with the norm
\begin{align}\label{def:HsT}
 \snorm{h}_{H^s(T)} := \inf\bigl\{ \snorm f_{H^s(\Rd)}: f\in H^s(\Rd) \mbox{ with } f_{|T} = h  \bigr\}\,.
\end{align}
Here, we followed \cite[Chapter 2.5.1]{EdTr96} and in particular we refer to the discussion there, when such restricted
spaces also  have an \emph{intrinsic} description.
In any case,
as inclusions between BSFs $E$ and $F$ on some set $T_0$ remain unchanged when both spaces are restricted to some subset $T\subset T_0$,
see e.g.~\cite[Proposition 3]{ScSt25a} for a general statement, the inclusions above lead to 
$H^s(T) \hookrightarrow H^r(T) \hookrightarrow  C_b(T)\hookrightarrow \sLx \infty T$
for all $s>r>d/2$. In this case, both $H^s(T)$ and $H^r(T)$ can thus be viewed to be  RKHSs.

With these preparations, we can now exactly describe those Gaussian processes that have an
RKHS that is essentially equivalent to a fractional Sobolev space and whose paths are contained in some RKHS.

\begin{theorem}\label{thm:sobol-main}
Let $T\subset \Rd$ be bounded and measurable with non-empty interior. Moreover, let 
$(\Om,\sA,\P)$ be a complete probability space and $X:= (X_t)_{t\in T}$ be a centered Gaussian process over  $(\Om,\sA,\P)$
whose  RKHS $H_X$ \changebegin consists of continuous functions and \changeend
is essentially equivalent to $H^s(T)$ for some $s>d/2$. Then the following statements are equivalent:
\begin{enumerate}
\item We have $s>d$.
\item There exist an RKHS  $H\subset \sLx \infty T$
and a version $Y$ of $X$ such that $\P(\{\sppath Y \in H\}) = 1$.
\changebegin 
\item There exist an RKHS  $H$ on $T$  
and a version $Y$ of $X$ with  $\P(\{\sppath Y \in H\}) = 1$.
\changeend
\end{enumerate}
In this case, for all $r \in (d/2, s-d/2)$ there actually exists 
such an $H$ that is essentially equivalent to  $H^r(T)$.
\end{theorem}

As the implication from \emph{ii)} to \emph{i)} is based upon  \changebegin Corollary \ref{cor:general-result} \changeend we see that the necessary condition 
derived there is sharp in the setup of Theorem  \ref{thm:sobol-main}. Below in Example \ref{ex:matern} 
we apply Theorem  \ref{thm:sobol-main} to the class of Mat\'ern kernels that are popular in many applications, see 
e.g.~\cite{RaWi06,PoBeScOa24a}.


Theorem \ref{thm:sobol-main} shows that we require a strong, dimension dependent regularity 
to ensure that the paths are contained 
in an RKHS. Interestingly, such a dimension dependent regularity disappears if we switch to (tensor) product kernels, 
which are used in e.g.~efficient implementations of Gaussian processes for machine learning \cite{GaPlWuWeWi18a, FaCoLoLiKiZh24a}.

%
%

\begin{theorem}\label{thm:mixed-sobol}
Let
$(\Om,\sA,\P)$ be a complete probability space, $T=(0,1)^d$, and $X:= (X_t)_{t\in T}$ be a centered Gaussian process over  $(\Om,\sA,\P)$ whose 
covariance function is of the form 
\begin{align*}
k_X(t,t') = \prod_{i=1}^d k_0(t_i,t_i')\, , \myqquad t,t'\in T,
\end{align*}
where $k_0:(0,1)\times (0,1)\to \R$ is a kernel with RKHS $H_0\smspac H^s(0,1)$ for some $s>1/2$. Then the following statements are equivalent:
\begin{enumerate}
\item We have $s>1$.
\item There exist an RKHS  $H\subset \sLx \infty T$
and a version $Y$ of $X$ such that $\P(\{\sppath Y \in H\}) = 1$.
\end{enumerate}
In this case, we may choose a Sobolev space of dominating mixed smoothness $H\cong B^r_{2,2}(0,1) \otimes_{2}\dots  \otimes_{2} B^r_{2,2}(0,1)$ for any  $r \in (1/2, s-1/2)$.
\end{theorem}

%% file: examples.tex
\section{Some Concrete Examples}\label{sec:examples}

In this section we present some explicit examples of Gaussian processes 
that do or do not have their paths in a suitable RKHS. 
To this end, 
we assume throughout this section that the considered Gaussian processes are defined  on some  complete
probability space $(\Om,\sA,\P)$. 
Moreover, $\lb$ and $\lbd$ denote the Lebesgue measures on the considered 
subsets $T\subset \R$, respectively $T\subset \R^d$.


With these preparations, our first example considers the Wiener process. 

\begin{example}\label{ex:wp}
Let $T:=[0,1]$ and
$W:= (W_t)_{t\in T}$ be the Wiener process. 
Then   $I_{W,\lb}:H_W\to \Lx 2 T$ is not nuclear and 
 \changebegin  there exists  no  RKHS  $H$ \changeend 
and  $\P(\{\sppath W \in H\}) = 1$. 
%
%
%
%
%
\end{example}

%
%

The proof of Example  \ref{ex:wp} shows that $I_{W,\lb}:H_W\to \Lx 2 T$ fails to be nuclear by a thin margin. 
This observation will be further investigated in Example \ref{ex:fbm}, which considers  the fractional Brownian motion.

Since the Brownian bridge $(B_t)_{t\in T}$ on $T:=[0,1]$ can be defined by $B_t := (1-t)W_t$, it is not surprising
that the negative result of Example \ref{ex:wp} carries over to the Brownian bridge. The details are formulated in the following example,
whose proof also illustrates the usefulness of the notion of essentially equivalent RKHSs.

\begin{example}\label{ex:bb}
Let $T:=[0,1]$ and $(B_t)_{t\in T}$ be the Brownian bridge.  
Then  \changebegin  there exists  no  RKHS  $H$ \changeend
such that  $\P(\{\sppath B \in H\}) = 1$. 
%
%
\end{example}

Like the Brownian bridge, the Ornstein-Uhlenbeck process is closely related to the Wiener process, and hence 
it is again not surprising that the negative result of Example \ref{ex:wp} carries over to 
the Ornstein-Uhlenbeck process. The following example presents the details for two variants of this process. 
%

\begin{example}\label{ex:ou-1}
Let $T:=[0,1]$ 
and 
$X:= (X_t)_{t\in T}$ be  a centered Gaussian process whose 
covariance function is given by either 
\begin{align*}
k_X^{(1)}(t_1,t_2) &:= a \exp(-\s |t_1-t_2|)\, ,  
\intertext{or}
k_X^{(2)}(t_1,t_2) &:= a \exp(-\s |t_1-t_2|) - a\exp(-\s(t_1+t_2)) \, ,  
\end{align*}
where  $a,\s>0$ are some constants.  Then 
\changebegin there exist no RKHS  $H$  \changeend
and no version $Y$ of $X$
 with $\P(\{\sppath Y \in H\}) =1$.
\end{example}


%

Our next example generalizes Example \ref{ex:wp} to the 
fractional Brownian motion. Here it turns out that the results depend on the Hurst parameter.

\begin{example}\label{ex:fbm}
Let $T:=[0,1]$ and $\a\in (0,1)$. Moreover, let 
 $B^{(\a)} = (B^{(\a)}_t)_{t\in T}$ be the fractional Brownian motion with Hurst index $\a$, that is, the centered Gaussian process with continuous paths and 
covariance function
\begin{align*}
k_{B^{(\a)}}(t_1,t_2) := \frac 12 \bigl( |s|^{2\a} + |t|^{2\a} - |t-s|^{2\a}  \bigr)  \, , \myqquad t_1,t_2\in [0,1]\, .
\end{align*}
Then there exists an 
\changebegin    RKHS  $H$ \changeend
with  $\P(\{\sppath B^{(\a)} \in H\}) = 1$ 
%
%
%
if and only if $\a>1/2$.
 Moreover, 
in the case $\a>1/2$ we may choose 
$H=H^s(T) \subset \sLx \infty T$ for any $s\in (1/2,\a)$. 
\end{example}


\changebegin
In the case $\a< 1/2$ the paths of the fractional Brownian motion are ``rougher''
than the paths of the Wiener process. Since we have already seen in Example 
\ref{ex:wp} that the paths of the Wiener process are not contained in an RKHS,
it is therefore not surprising that the paths of $B^{(\a)}$ are not contained in an 
RKHS, either. 
In the case $\a>1/2$, it is well-known that the paths of  $B^{(\a)}$ are $\b$-H\"older-continuous for all $\b\in (0,\a)$. The positive result then follows by standard 
embeddings between Besov spaces, see the proof of Example \ref{ex:fbm} for details.
\changeend

The fractional Brownian motion is closely related to the Riemann-Liouville process, see e.g.~\cite[Example 3.4]{Lifshits12} or \cite{Picard11a}, which is considered in the following example. 

\begin{example}\label{ex:rlp}
Let $T:=[0,1]$ and $\a\in (0,1)$.
Moreover, let $R^{(\a)} = (R^{(\a)}_t)_{t\in T}$  be
%
the Riemann-Liouville process with index $\a$, that is, the centered Gaussian process with continuous paths given by the It\^o-Integral
\begin{align*}
R^{(\a)}_t := \frac 1{\Gamma{(\a+1/2})}  \int_0^t (t-s)^{\a-1/2} \intd W_s\, , \myqquad t\in [0,1].
\end{align*}
Then there exists an 
\changebegin    RKHS  $H$  \changeend
with  $\P(\{\sppath R^{(\a)} \in H\}) = 1$ 
%
%
if and only if $\a>1/2$.
Moreover, in the case $\a>1/2$ we may choose $H=H^s(T)\subset \sLx \infty T$ for any $s\in (1/2,\a)$.
\end{example}

The next example, which  considers Gaussian processes on a higher dimensional $T$, uses  our result for RKHSs of Sobolev type 
formulated in Theorem \ref{thm:sobol-main}. Analogous results can thus be obtained for other kernels of Sobolev type,
such as the Wendland kernels, which are often attractive from a computational point of view, see \cite[Chapter 9]{Wendland05}.

\begin{example}\label{ex:matern}
    Let $T\subset \Rd$  be  bounded and measurable with non-empty interior.
%
Moreover, let $X:= (X_t)_{t\in T}$ be a centered Gaussian
process whose  covariance function is a Mat\'ern kernel of order $\a>0$ and width $\s>0$, that is 
\begin{align*}
 k_X(t_1,t_2) = \bigl( \s \snorm{t_1-t_2} \bigr)^\a h_\a\bigl( \s \snorm{t_1-t_2} \bigr)\, , \myqquad t_1,t_2\in T,
\end{align*}
where $h_\a$ denotes the modified Bessel function of the second type of order $\a$. Then there 
exist 
\changebegin  an  RKHS  $H$  \changeend
and a version $Y$ of $X$ with  $\P(\{\sppath Y \in H\}) =1$
if and only if $\a>d/2$.
Moreover, in the case $\a>d/2$ we may choose   $H=H^r(T)$ for any $r\in (d/2,\a)$.
\end{example}

\changebegin
To compare the smoothness of the paths   and the elements of the RKHS 
in Example \ref{ex:matern}, we quickly assume that $\a>d/2$. Then 
the paths of $X$ are contained in $H^{\a-\e}(T)$ for all sufficiently small $\e>0$.
Moreover, 
the proof of Example \ref{ex:matern} uses that the 
 RKHS $H_{\a,\s}$ of the Mat\'ern kernel equals $H^{\a + d/2}$ up to equivalent norms.
Consequently, there is a $d/2+\e$-gap, when the smoothness is measured in the Sobolev sense. This gap coincides with the one observed in 
\cite{Steinwart19a}, where it is also shown that it is essentially sharp, 
see also \cite{Karvonen23a} for a more detailed analysis and \cite{Scheuerer10a} for 
similar results for weakly stationary processes.  Finally, \cite[Example 4.8]{Steinwart19a} shows that in the case $d=1$ these Sobolev results 
improve classical results 
on the continuous  differentiability of the paths of $X$, of e.g.~\cite{CrLe67}.
\changeend

Our last example considers stationary processes on compact Abelian groups $T$ with Haar measure $\lb$. 
It mostly relies on the fact that for translation invariant, $\Cfield$-valued kernels the characters of $T$
form on ONB of $\Lx 2 {\lb,\Cfield}$ consisting of eigenvectors of $T_{k,\lb}$.

\begin{example}\label{ex:abel}
 Let $(T,+)$ be a compact Abelian group with Haar measure $\lb$
 and $X:= (X_t)_{t\in T}$ be a  centered Gaussian process whose covariance 
 function $k_X$ is measurable and satisfies
 \begin{align*}
  k_X(t_1,t_2) = k_X(t_1-t_2,0)\, , \myqquad t_1,t_2\in T.
 \end{align*}
Then there exist an RKHS $H\subset \sLx \infty \sB$ and 
a version $Y$ of $X$  with $\P(\{\sppath Y \in H\}) =1$ if and only if 
\begin{align*} 
\sum_{i\in I} \sqrt{\mu_i(T_{X,\lb})} < \infty.
\end{align*}
Moreover, in this case we may choose $H=H_X^{1/2}$ and we have $H_X^{1/2} \subset C(T)$.
\end{example}

Since Example \ref{ex:abel}  is an application of Theorem \ref{thm:ew-char}, it 
 can be easily extended to situations in which the spectral properties of $k_X$ are well understood
 in the sense of the discussion following Theorem \ref{thm:ew-char}. 
 Such situations 
  include
 e.g.~isotropic Gaussian random fields on 
$\mathbb S^2$, and with some more work also on $\mathbb S^d$,
if $\lb$ is the surface measure.  We refer to \cite{LaSc15a} for the details on the required spectral properties \changebegin and to \cite{HoSc26a} for an application for neural networks.\changeend

%% file: proofs.tex
\section{Proofs}\label{sec:proofs}

 \subsection{Hilbert-Schmidt and Nuclear Operators}\label{subsec:hs+nuc-op}
 
 In this subsection we briefly recall  some aspects of compact operators between Hilbert spaces. 
To this end, let $H$ be a Hilbert space and $R:H\to H$ be a  compact, positive, and self-adjoint operator. 
Then the  spectral theorem for such operators,  e.g.~\cite[Chapter V.2.3]{Kato80} or \cite[Theorem VI.3.2]{Werner11}, shows that
we have at most countably many non-zero eigenvalues $(\mu_i(R))_{i\in I}$, where we have set either $I:=\{1,\dots,n\}$ or $I:=\N$, 
and where we have included the geometric multiplicities. Moreover, there exists an ONS $(e_i)_{i\in I}$ in $H$ such that
$e_i$ is an eigenvector 
of $R$ to the eigenvalue $\mu_i(R)$ for all $i\in I$.
Finally, these eigenvalues are positive and in the case $I=\N$ they can be ordered to  form a decreasing sequence converging to $0$.
In the following, we always assume a decreasing order. \changebegin 
Finally, the spectral theorem leads to the following representation 
\begin{align}\label{eq:spectral-formula}
Rx = \sum_{i\in I} \mu_i(R) \skprod x{e_i} e_i \, , \myqquad x\in H.
\end{align}
In particular, we have $H = \ker R \oplus    \overline{\spann\{e_i: i\in I\}}^{H}$.
\changeend

If $S:H_1\to H_2$ is a compact operator between two Hilbert spaces $H_1$ and $H_2$, we write 
$S^*:H_2\to H_1$ for the adjoint of $S$. Then 
the operator $S\circ S^*:H_2\to H_2$
is compact, positive, and self-adjoint, and hence we can apply the spectral theorem to $S\circ S^*$. 
By doing so, we can define the
singular numbers of $S$   by
\begin{align*}
s_i(S) := \sqrt{\mu_i(S\circ S^*)} \, , \myqquad i\in I,
\end{align*}
and $s_i(S) := 0$ if $I:=\{1,\dots,n\}$ and $i>n$.
Recall that for compact, positive, and self-adjoint $R:H\to H$
we have 
\begin{align}\label{eq:sing-vs-eigen}
s_i(R) = \sqrt{\mu_i(R\circ R^*)}  = \sqrt{\mu_i(R^2)} =    \mu_i(R) \, ,
\end{align}
 where the last identity directly follows from the spectral series representation, see e.g.~\cite[Equation (2.21)  in Chapter V.2.3]{Kato80} or  \cite[Theorem VI.3.2]{Werner11}.
For compact $S:H_1\to H_2$ this implies
\begin{align}\label{eq:snum-factor}
s_i^2(S) = \mu_i(S\circ S^*) = s_i(S\circ S^*)\, ,
\end{align}
and in addition, there exist  ONSs $(e_i)_{i\in I}$ in $H_1$ and $(f_i)_{i\in I}$ in $H_2$ with 
\begin{align}\label{eq:schmidt-repr}
Sx = \sum_{i\in I} s_i(S) \skprod x{e_i}_{H_1} f_i \, ,\myqquad  x\in H_1,
\end{align}
see e.g.~\cite[Equation (2.23) in  Chapter V.2.3]{Kato80} or \cite[Theorem VI.3.6]{Werner11}.
Finally, recall that a compact $S:H_1\to H_2$ is nuclear if and only if
we have 
\begin{align}\label{eq:nuc-op}
\sum_{i\in I} s_i(S) < \infty\,   ,
\end{align}
see e.g.~\cite[Chapter 11.2]{BiSo87}. Moreover, by definition $S$ is a Hilbert-Schmidt operator
if $\sum_{j\in J} \snorm{Se_j}_{H_2}^2 < \infty$ for some (and then all)
ONB(s) $(e_j)_{j\in J}$ of $H_1$, or equivalently, 
if and only if 
\begin{align}\label{eq:hs-op}
\sum_{i\in I} s_i^2(S) < \infty\, ,  
\end{align}
see e.g.~\cite[Chapter 11.3, Theorem 1]{BiSo87}. By \eqref{eq:snum-factor} we then immediately obtain the following 
elementary yet useful lemma.

\begin{lemma}\label{lem:hs-vs-nuc}
Let $S:H_1\to H_2$ be a compact operator between two Hilbert spaces $H_1$ and $H_2$. 
Then $S$ is a 
Hilbert-Schmidt operator if and only if $S\circ S^*$ is nuclear.
\end{lemma}

 Finally, it is sometimes helpful to compare the singular numbers with some other $s$-number scales in the sense of \cite[Chapter 2]{Pietsch87}.
 To this end, let us first recall that for a given bounded linear operator $S:E\to F$ between two Banach spaces $E$ and $F$,
 the $n$-th  approximation number is defined by 
 \begin{align}\label{eq:an}
 a_n(S):= \inf\bigl\{\snorm {S-A}\,\bigl|\, A:E\to F \mbox{ bounded and linear  with } \rank A<n\bigr\} \, .
 \end{align}
 Furthermore,
the $n$-th Kolmogorov number is defined by 
 \begin{align}\label{eq:dn}
d_n(S) := \inf\bigl\{\e>0:  \exists F_\e \mbox{ subspace of $F$ with $\dim F_\e < n$ and }  SB_E \subset F_\e+\e B_F   \bigr\}\,,
\end{align}
where, $B_E$ and $B_F$ denote the closed unit balls in $E$ and $F$ and we followed the exposition in \cite[Chapter 2.2]{CaSt90}.
It is well known that  the approximation
and 
 Kolmogorov numbers form an $s$-scale
in the sense of \cite[Chapter 2.2.1]{Pietsch87}. We refer to \cite[Chapters 2.3 and 2.5]{Pietsch87} and \cite[Chapters 2.1 and 2.2]{CaSt90} for corresponding results. In addition, the 
Kolmogorov numbers are a surjective $s$-scale, see e.g.~\cite[Chapter 2.5]{Pietsch87}. 
Moreover, $S:E\to F$ is compact if and only if  $d_n(S) \to 0$, see 
e.g.~\cite[page 49]{CaSt90}. Furthermore, for operators $S:H_1\to H_2$   between Hilbert spaces
we have $a_n(S) =  d_n(S)$, see e.g.~\cite[Theorem 11.3.4]{Pietsch80}, and if  
$S:H_1\to H_2$ is compact we even have 
\begin{align}\label{eq:kol-vsr-sing}
s_n(S) = a_n(S) =   d_n(S)\, , \myqquad n\geq 1,
\end{align}
by \cite[Proposition 11.3.3]{Pietsch80} as \eqref{eq:schmidt-repr} is a so-called Schmidt representation in 
the sense of \cite[D.3.2]{Pietsch80}.
%
In other words, we may freely switch between these  quantities as soon as we know that $S$ is compact.


 \subsection{Nuclear Dominance}

 In this subsection we show that the notion of nuclear dominance used in \cite{LuBe01a} 
 is equivalent to ours. To this end, we begin with recalling \cite[Theorem 1.1]{LuBe01a}.

 \begin{theorem}\label{thm:dominance}
Let $H_1$ and $H_2$ be RKHSs over $T$ with  $H_1\subset H_2$.
Then the resulting embedding map
\begin{align}\label{eq:embed-map}
J:H_1&\to H_2 \\ \nonumber
h_1&\mapsto h_1
\end{align}
is linear and bounded. Moreover, there exists a unique linear $D_{H_1,H_2}:H_2\to H_2$ that satisfies both
 $\ran D_{H_1,H_2}\subset H_1$
and
\begin{align}\label{thm:dominance-h1}
\skprod {h_2}{h_1}_{H_2} = \skprod {D_{H_1,H_2} h_2}{h_1}_{H_1}\, , \myqquad h_1\in H_1,h_2\in H_2.
\end{align}
In addition, $D_{H_1,H_2}$ is bounded, positive and self-adjoint.
\end{theorem}

Recall that \cite{LuBe01a} speaks of nuclear dominance if $D_{H_1,H_2}$ is nuclear. 
The following simple lemma shows that this notion is equivalent to our Definition \ref{def:nuc-dom}.

\begin{lemma}\label{lem:comput-dominance}
Let $H_1$ and $H_2$ be RKHSs over $T$ with  $H_1\subset H_2$ and $J:H_1\to H_2$ be
the corresponding embedding map \eqref{eq:embed-map}.
Then we have 
\begin{align*}
D_{H_1,H_2} = J\circ J^*\, .
\end{align*}
In particular, $D_{H_1,H_2}$ is nuclear if and only if $J$ is a Hilbert-Schmidt operator.
\changebegin If, in this case, $(\mu_i)_{i\in I}$ are the at most countably many non-zero eigenvalues
(including geometric multiplicities) of $D_{H_1,H_2}$ and $(f_i)_{i\in I}\subset H_2$ is an ONS of 
corresponding eigenvectors, then we have 
\begin{align}\label{eq:spectral-orth-sum}
H_2 = \ker D_{H_1,H_2} \oplus \overline{\spann\{f_i: i\in I\}}^{H_2} = \ker D_{H_1,H_2} \oplus \overline {\ran J}^{H_2}\, .
\end{align}
Moreover, the vectors $e_i := \mu_i^{-1} J^*f_i\in H_1$ satisfy $Je_i = f_i$ and $(\sqrt{\mu_i} e_i)_{i\in I}$ is an ONS in $H_1$. 
\changeend
\end{lemma}

\begin{proofof}{Lemma \ref{lem:comput-dominance}}
Let us write  $R:= J \circ J^* :H_2\to H_2$.
By construction, its range satisfies $\ran R\subset H_1$. Moreover, for $h_1\in H_1$ and $h_2\in H_2$ we have
\begin{align*}
\skprod {h_2}{h_1}_{H_2} = \skprod {h_2}{J h_1}_{H_2} = \skprod {J^*h_2}{h_1}_{H_1} = \skprod {(J \circ J^*)h_2}{h_1}_{H_1}
= \skprod {R h_2}{h_1}_{H_1} \, .
\end{align*}
Consequently, $R$ satisfies \eqref{thm:dominance-h1}, and  by uniqueness stated in  Theorem \ref{thm:dominance} we can thus conclude $R=D_{H_1,H_2}$.
The \changebegin equivalence \changeend is a direct consequence of Lemma \ref{lem:hs-vs-nuc} and $D_{H_1,H_2} = J\circ J^*$.

\changebegin
Moreover, the first identity in \eqref{eq:spectral-orth-sum} holds as noted below
 \eqref{eq:spectral-formula}. To establish the second identity
we first show
\begin{align}\label{lem:comput-dominance-h1}
\ker J^* = \ker D_{H_1,H_2}\, .
\end{align}
Here the inclusion ``$\subset$'' follows from $D_{H_1,H_2} = J\circ J^*$, and  
the converse inclusion follows from 
\begin{align*}
\snorm{J^*h_2}_{H_1} = \skprod {J^*h_2}{J^*h_2}_{H_1}  = \skprod {D_{H_1,H_2}h_2}{h_2}_{H_2}  = 0\, , \myqquad h_2\in \ker D_{H_1,H_2}.
\end{align*}
Now, combining the first identity in  \eqref{eq:spectral-orth-sum} with \eqref{lem:comput-dominance-h1} yields
\begin{align}\label{lem:comput-dominance-h2}
\overline{ \spann\{f_i: i\in I\}}^{H_2} = (\ker D_{H_1,H_2})^\perp = (\ker J^*)^\perp = \overline{\ran J}^{H_2}\, , 
\end{align}
where the last equation is a textbook formula, see e.g.~\cite[Theorem II.2.19 with Corollary I.2.10]{Conway90}.

Finally,  the definition of the $e_i$ immediately yields 
\begin{align*}
Je_i = \mu_i^{-1} (J\circ J^*)f_i = \mu_i^{-1} D_{H_1,H_2} f_i = f_i\, , 
\end{align*}
as well as 
\begin{align*}
\langle e_i, e_j\rangle_{H_1} 
= \mu_i^{-1} \mu_j^{-1}    \langle  J^*f_i , J^*f_j\rangle_{H_1} 
= \mu_i^{-1} \mu_j^{-1}    \langle  f_i , D_{H_1,H_2} f_j\rangle_{H_2} 
= \mu_i^{-1}      \langle  f_i ,   f_j\rangle_{H_2} \, .
\end{align*}
With the help of this identity we easily see that $(\sqrt{\mu_i} e_i)_{i\in I}$ is an ONS in $H_1$. 
%
%
\changeend
\end{proofof}


\subsection{Choosing Small RKHSs}

\changebegin

In this subsection, we show how we can ``shrink'' the RKHS containing the paths 
of a suitable version of a Gaussian process in the situation of Theorem \ref{thm:LuBe01a}.

To this end, 
let us recall that a map $\xi:\Om\to H$  from a measurable space $(\Om,\sA)$ 
into a Hilbert space $H$ is weakly measurable, if $\skprod \xi f_H:\Om\to \R$
is measurable for all $f\in H$. Moreover, it is called strongly measurable
if it is weakly measurable and $\xi(\Om)$ is separable in $H$. In this case,
the function $\om\to \snorm {\xi(\om)}_H$ is also measurable.

Now, our first result shows that we can always assume that the space $H$ is separable as
long as we are willing to potentially switch the versions of $X$.

\begin{lemma}\label{lem:lukic-beder-covariance}
Let $(\Om,\sA,\P)$ be a complete probability space and $X:= (X_t)_{t\in T}$ be a centered Gaussian process over  $(\Om,\sA,\P)$
with  RKHS $H_X$. Moreover, let 
$H$ be an RKHS on $T$  with kernel $k$ and $Y$ be a version of $X$  with
\begin{align*}
\P(\{\sppath Y \in H\}) >0\, .
\end{align*}
Then there exists a separable RKHS $H^\dagger \subset H$ and another version $Y^\dagger$ of $X$
with $\sppath Y^\dagger(\om) \in H^\dagger$ for all $\om \in \Om$ 
such that the resulting map $\Om\to H^\dagger$ given by $\om\mapsto \sppath Y^\dagger(\om)$
is strongly measurable and satisfies 
\begin{align*}
\E_\P |\langle \sppath Y^\dagger, f\rangle_{H^\dagger}|^2 = \langle D_{H_X,H^\dagger}f, f\rangle_H\,, \myqquad f\in H^\dagger. 
\end{align*}
\end{lemma}

\begin{proofof}{Lemma \ref{lem:lukic-beder-covariance}}
 Theorem \ref{thm:LuBe01a} gives  $H_X\ll H$ and a version $\tilde Y$ of $X$ with 
 $\P(\{\sppath {\tilde Y} \in H\}) =1$. Obviously, we may assume without loss of generality that 
 this version actually satisfies $\sppath {\tilde Y}(\om) \in H$ for all $\om \in \Om$.
 
Let us now define $\xi:\Om\to H$  
by $\xi(\om) := \sppath {\tilde Y}(\om)$ for all $\om \in \Om$. Then 
\cite[Lemma 2.1]{LuBe01a} shows that $\xi$ is weakly measurable and 
the reproducing property of $k$ yields
\begin{align}\label{lem:lukic-beder-covariance-h1}
\tilde Y_t (\om) = \langle \xi(\om), k(t,\mycdot)\rangle_{H}\, , \myqquad \om\in \Om, t\in T.
\end{align}
 Moreover, since $\tilde Y$ is a version of $X$ we have $H_{\tilde Y} = H_X\ll H$. 
Consequently
\cite[Corollary 3.1]{LuBe01a}
gives a strongly measurable $\xi^\dagger:\Om\to H$ with $\snorm {\xi^\dagger}_H \in \sLx 2 \P$
and 
\begin{align}\label{lem:lukic-beder-covariance-h2}
\P\bigl(\bigl\{ \om\in \Om :  \langle \xi^\dagger(\om), f\rangle_{H} = \langle \xi(\om), f\rangle_{H}   \bigr\}   \bigr) = 1 \, , \myqquad f\in H.
\end{align}
Alternatively, one can obtain $\xi^\dagger$ by \cite[Theorem 5.1]{Edgar77a}, see also \cite[page 89]{DiUh77} for a more detailed proof, 
as $\xi$ is Pettis integrable and Hilbert spaces have the Radon-Nikodym property.

We now define $H^\dagger := \overline{\xi^\dagger(\Om)}^{H}$. Since $\xi^\dagger$ is strongly measurable, $H^\dagger$ is separable. In addition, as a closed subspace of an RKHS it is 
an RKHS. Let us now consider 
 process $Y^\dagger := (Y^\dagger_t)_{t\in T}$ defined by 
\begin{align*}
Y^\dagger_t (\om) := \langle \xi^\dagger(\om), k (t,\mycdot)\rangle_{H }\,, \myqquad t\in T,\om\in \Om.
\end{align*}
Then our construction gives 
$\sppath Y^\dagger(\om) = \xi^\dagger(\om) \in H^\dagger$ for all $\om \in \Om$.
Moreover, combining this definition with \eqref{lem:lukic-beder-covariance-h2}
and \eqref{lem:lukic-beder-covariance-h1} we see that $Y^\dagger$ is a version of $\tilde Y$, and thus also a version of $X$.

Now, since $Y^\dagger$ is a version of $X$ with $\P(\{\sppath Y^\dagger \in H^\dagger\}) =1$
another application of Theorem \ref{thm:LuBe01a} yields
$H_{Y^\dagger} = H_X \ll H^\dagger$. 
Applying 
 \cite[Corollary 3.1]{LuBe01a}, see also  \cite[Equation (3.2)]{LuBe01a}, to $Y^\dagger$
 then gives the last assertion.
\end{proofof}

The next result, which is the main theorem of this subsection, refines the
``shrinking'' to a somewhat minimal space.

\begin{theorem}\label{thm:nice-rkhs}
Let $(\Om,\sA,\P)$ be a complete probability space and $X:= (X_t)_{t\in T}$ be a centered Gaussian process over  $(\Om,\sA,\P)$
with  RKHS $H_X$. Moreover, let 
$H$ be an RKHS on $T$ and $Y$ be a version of $X$  with  
\begin{align*}
\P(\{\sppath Y \in H\}) >0\, .
\end{align*}
Then there exist a version $Y^*$ of $X$ and a  separable RKHS $H^*$ on $T$ such that $H_X\subset H^*$ is dense in $H^*$  and 
\begin{align*}
\P(\{\sppath Y^* \in H^*\}) =1 \, .
\end{align*}
\end{theorem}

\begin{proofof}{Theorem \ref{thm:nice-rkhs}}
Without loss of generality we may assume that $Y=Y^\dagger$ and $H=H^\dagger$, where 
$Y^\dagger$ and $H^\dagger$ are according to 
Lemma \ref{lem:lukic-beder-covariance}. In particular, we may assume that $H$ is separable.
In the following we consider $Y^* := Y$, so that it remains to find $H^*$.

Since  Theorem \ref{thm:LuBe01a} ensures $H_X\ll H$, 
Lemma \ref{lem:comput-dominance} gives us 
an at most countable ONS $(f_i)_{i\in I}\subset H$ 
consisting of eigenvectors of $ D_{H_X, H} $ with non-zero eigenvalues $(\mu_i)_{i\in I}$
such that 
\begin{align}\label{thm:nice-rkhs-h1}
H = \ker D_{H_X, H} \oplus  \overline{\spann\{f_i: i\in I\}}^{H}\, .
\end{align}
We define $H^* :=  \overline{\spann\{f_i: i\in I\}}^{H}$ and equip it with the norm of $H$. 
Then $H^*$ is a separable RKHS  and 
the second identity of \eqref{eq:spectral-orth-sum} in 
Lemma \ref{lem:comput-dominance} shows 
\begin{align*}
H^* =  \overline{\spann\{f_i: i\in I\}}^{H} = \overline {\ran J}^H = \overline {H_X}^H  \, , 
\end{align*}
where the third identity  is true since $J:H_X\to H$ is the embedding map. In particular, 
we have $H_X\subset H^*$, and for every $h\in H^*$ there exists a sequence 
$(f_n)\subset H_X$ with $\snorm{f_n-h}_H \to 0$. Since we know $f_n\in H^*$ we conclude that 
$\snorm{f_n-h}_{H^*} \to 0$. In other words, we have shown $H^* \subset \overline {H_X}^{H^*}$.

Now, since $H$ is separable, so is its subspace $\ker  D_{H_X,H}$. Consequently, there
 exists a countable dense subset $\frk D\subset \ker  D_{H_X,H}$.
For $f\in \frk D$ the last assertion of Lemma   \ref{lem:lukic-beder-covariance} yields
\begin{align*}
\E_\P |\langle \sppath Y, f\rangle_{H}|^2 = \langle D_{H_X,H}f, f\rangle_H = 0\, , 
\end{align*}
and hence there exists an $N_f\in \sA$ with $\P(N_f) = 0$ and  $\langle \sppath Y(\om), f\rangle_H = 0$ for all $\om \in \Om\setminus N_f$. Let us define 
\begin{align*}
N := \bigcup_{f\in \frk D} N_f\, .
\end{align*}
Clearly, we have $N\in \sA$ with $\P(N) = 0$ and $\langle \sppath Y(\om), f\rangle_H = 0$ for all $f\in \frk D$ and $\om \in \Om\setminus N$.

With these preparations we can now show $\sppath Y(\om) \in H^*$ for all $\om \in \Om\setminus N$. To this end, we pick an $f\in \ker D_{H_X,H}$. Then there exists a sequence 
$(f_n)\subset \frk D$ with $f_n\to f$ in $H$. For 
$\om \in \Om\setminus N$ our construction then shows
\begin{align*}
\langle \sppath Y(\om), f\rangle_H
= \lim_{n\to \infty} \langle \sppath Y(\om), f_n\rangle_H = 0\, .
\end{align*}
In other words, we have found $\sppath Y(\om) \in ( \ker D_{H_X,H})^\perp = H^*$,
where the last identity follows from \eqref{thm:nice-rkhs-h1} and the definition of $H^*$.
\end{proofof}

\changeend


\subsection{Proof of Theorem \ref{thm:general-result-envelope} and its Corollaries}

Before we can prove 
\changebegin Theorem \ref{thm:general-result-envelope} and its corollaries \changeend
we need to establish \changebegin a few auxiliary results.
We begin with an elementary lemma describing the effect of multiplying a kernel by some function.

\begin{lemma}\label{lem:mult-of-rkhs}
Let $k$ be a kernel on $T$ with RKHS $H$ and $\a:T\to \R$ be a function. Then 
the RKHS of the kernel
$\a \star k:T\times T\to \R$ defined by 
\begin{align*}
(\a \star k)(t,t') := \a(t) \a(t') k(t,t')\, , \myqquad t,t'\in X,
\end{align*}
is $H_{\a \star k} := \{ \a h: h\in H  \}$ with norm given by $\snorm{f}_{H_{\a \star k}} = \inf\bigl\{ \snorm h_H:  h\in H \mbox{ with } f = \a h  \bigr\}$.
\end{lemma}

\begin{proofof}{Lemma \ref{lem:mult-of-rkhs}}
Let $\Phi:T\to H$ be the canonical feature map of $k$, that is $\Phi(t) = k(\mycdot, t)$ for all $t\in T$.
For $t_1,t_2\in T$ an elementary calculation then shows 
\begin{align*}
\skprodb{\a(t_1)\Phi(t_1)}{\a(t_2)\Phi(t_2)}_H
\scriptonlyraw{&}= \a(t_1) \a(t_2) \skprod{\Phi(t_1)}{\Phi(t_2)}_H  
\scriptonlyraw{&}=  \a(t_1) \a(t_2) k(t_1,t_2)\, .
\end{align*}
Consequently, $\a \Phi:T\to H$ is a feature map of $\a \star k$. Moreover,   for $h\in H$ and $t\in T$
another simple calculation shows
\begin{align*}
\skprodb h {(\a \Phi)(t)}_H = \a(t) \, \skprod h {\Phi(t)}_H = \a(t) h(t)\, .
\end{align*}
Now the assertion follows  from the fundamental theorem of RKHSs, see e.g.~\cite[Theorem 4.21]{StCh08}.
\end{proofof}

The next auxiliary results describes the relationship between envelope functions and the behavior of the kernel on its diagonal.

\begin{lemma}\label{lem:envelope-vs-kernel}
Let $H$ be an RKHS on $T$ with kernel $k$. Then   $b:T\to [0,\infty)$  is an envelope for $H$, if and only if 
\begin{align}\label{lem:envelope-vs-kernel-eq}
\sqrt{k(t,t)} \leq b(t) \, , \myqquad t\in T.
\end{align}
\end{lemma}

\begin{proofof}{Lemma \ref{lem:envelope-vs-kernel}}
If \eqref{lem:envelope-vs-kernel-eq} is true, then the assertion immediately follows from \eqref{eq:all-envelope}.

To establish the converse implication, we define $\a :T\to [0,\infty)$ by 
\begin{align*}
\a(t) := 
\begin{cases}
1/b(t), & \mbox{ if } b(t) > 0 \, ,\\
0 , & \mbox{ otherwise.}
\end{cases} 
\end{align*}
Let us now fix an $f\in H_{\a \star k}$, where $H_{\a \star k}$ is the RKHS considered in 
 Lemma \ref{lem:mult-of-rkhs}.
According to Lemma \ref{lem:mult-of-rkhs} there then exists an $h\in H$ with $f=\a h$. 

Our next, intermediate goal is to show that there is only one such  $h$. To this end, we note that for $t\in T$ with $b(t) >0$ we obviously have 
$h(t) = b(t) f(t)$ by the definition of $\a$. Moreover, for 
$t\in T$ with $b(t) = 0$ the assumed envelope condition \eqref{eq:envelope} ensures   $h(t) = 0$. 
In other words, the function $h$ is uniquely determined for all $t\in T$, and therefore Lemma 
 \ref{lem:mult-of-rkhs} shows 
 \begin{align*}
  \snorm{f}_{H_{\a \star k}} =  \snorm h_H \, .
 \end{align*}

With these preparations we now observe that, for all $t\in T$, the envelope condition gives 
\begin{align*}
|f(t)| = |\a(t)| \cdot |h(t)| \leq    |\a(t)| \cdot \snorm h_H\cdot   b(t) \leq  \snorm h_H = \snorm{f}_{H_{\a \star k}}\, .
\end{align*}
Since $f\in H_{\a \star k}$ was arbitrary, we have thus found $H_{\a \star k} \hookrightarrow \ell_\infty(T)$ with 
 $\snorm{\id: H_{\a \star k} \to \ell_\infty(T)} \leq 1$. By standard RKHS theory, see e.g.~\cite[Lemma 4.23]{StCh08}, we then obtain 
 \begin{align*}
 \a(t) \a(t) k(t,t) = ( \a \star k)(t,t) \leq 1 \, , \myqquad t\in T.
 \end{align*}
 For $t\in T$ with $b(t) >0$ this obviously implies \eqref{lem:envelope-vs-kernel-eq}. Moreover, for $t\in T$ with $b(t) = 0$ the envelope condition 
 applied to $h := k(\mycdot ,t)\in H$ gives 
 \begin{align*}
 |h(t')| \leq  \snorm h_H \cdot b(t')  \, , \myqquad t'\in T, 
 \end{align*}
 and for $t':= t$ the latter implies $k(t,t) = h(t) = 0$. Consequently,  \eqref{lem:envelope-vs-kernel-eq} is also true for $t\in T$ with $b(t) = 0$.
\end{proofof}

Our next and final auxiliary result investigates the operator $I_{k,\nu}:H\to \Lx 2\nu$ in more detail.

\changeend

\begin{lemma}\label{lem:Iknu-new}
Let   $k:T\times T\to \R$ be a kernel with RKHS $H$ and  $\sB$ be a $\s$-algebra on $T$ \changebegin with  $H\subset \sLx 0\sB$.
Moreover, let $b:T\to [0,\infty)$ be a measurable envelope function for $H$. Then for every measure $\nu$ on $(T,\sB)$
with $b\in \sLx 2 \nu$ we have $H\subset \sLx 2 \nu$ and the resulting \changeend
``embedding'' operator
\begin{align}\label{eq:H2L2-new}
I_{k,\nu}:H&\to \Lx 2 \nu \\ \nonumber
h &\mapsto \eqclass h  
\end{align}
is  linear and Hilbert-Schmidt. Moreover, its adjoint $I^*_{k,\nu}:\Lx 2 \nu\to H$ is given by 
\begin{align*}
(I_{k,\nu}^* f)(s) =  \int_T k(s, t) f(t) \intd \nu(t)\, , \myqquad f\in \Lx 2 \nu,\, s\in T.
\end{align*}
\end{lemma}

\begin{proofof}{Lemma \ref{lem:Iknu-new}}
\changebegin We first note that the envelope condition \eqref{eq:envelope}  gives $|h(t)|^2 \leq \snorm h_H^2 \cdot b^2(t)$ for all $h\in H$ and $t\in T$.
By the assumed $b\in \sLx 2 \nu$ we thus find $H\subset \sLx 2 \nu$. This in turn shows that $I_{k,\nu}$ is well-defined. Moreover, the linearity of 
$I_{k,\nu}$ is obvious.

To verify that $I_{k,\nu}$ is Hilbert-Schmidt, we fix an ONB $(e_j)_{j\in J}$ of $H$. For finite $A\subset J$ we then find
\begin{align*}
\sum_{j\in A} \snorm {I_{k,\nu} e_j}_{\Lx 2 \nu}^2 
= 
\int_T \sum_{j\in A} e^2_j(t) \intd \nu(t)
= 
\int_T \sum_{j\in A} \bigl| \skprod {e_j}{k(\mycdot, t)}_H\bigr|^2 \intd \nu(t) \, .
\end{align*}
Moreover, Bessel's inequality gives 
\begin{align*}
 \sum_{j\in A} \bigl| \skprod {e_j}{k(\mycdot, t)}_H\bigr|^2 
 \leq \snorm{k(\mycdot, t)}_H^2 
 = k(t,t)
 \leq b^2(t)\, ,
\end{align*}
where in the last step we used Lemma \ref{lem:envelope-vs-kernel}. Combining this estimate with our initial identity then yields
\begin{align*}
\sum_{j\in A} \snorm {I_{k,\nu} e_j}_{\Lx 2 \nu}^2 
= 
\int_T \sum_{j\in A} \bigl| \skprod {e_j}{k(\mycdot, t)}_H\bigr|^2 \intd \nu(t)
\leq 
\snorm b_{\sLx 2 \nu}^2\, .
\end{align*}
Since this is true for all finite $A\subset J$ it is also true for $A:= J$, that is, $I_{k,\nu}$ is Hilbert-Schmidt.

 \changeend

Let us finally derive the formula for the adjoint operator. To this end, we fix an $f\in \Lx 2 \nu$ and an $s\in T$. We then obtain 
\begin{align*}
(I_{k,\nu}^* f)(s) 
= \skprodb{I_{k,\nu}^* f}{k(s,\mycdot)}_H 
= \skprodb{f}{I_{k,\nu} (k(s,\mycdot))}_{\Lx 2 \nu}
= \skprodb{f}{ \eqclass{k(s,\mycdot)}}_{\Lx 2 \nu}
= \int_T k(s, t) f(t) \intd \nu(t)\, ,
\end{align*}
where in the first step we used the reproducing property $k$.
\end{proofof}

\begin{proofof}{Theorem \ref{thm:general-result-envelope}}
\changebegin  By Theorem \ref{thm:nice-rkhs} we may assume without loss of generality that 
$H_X\subset H$ is dense in $H$. Since convergence in $H$ implies pointwise
convergence and we assumed $H_X\subset \sLx 0 \sB$, we therefore also find 
$H\subset \sLx 0 \sB$. In turn,
$b\in \sLx 2 \nu$  in combination with Lemma  \ref{lem:Iknu-new} ensures that  $I_{k,\nu}$ 
is   a Hilbert-Schmidt operator, where 
$k$ is the kernel of $H$ and $I_{k,\nu}$ is defined by \eqref{eq:H2L2-new}.

Moreover, Theorem \ref{thm:LuBe01a} ensures 
$H_X\ll H$, that is the embedding map $J:H_X\to H$ is a Hilbert-Schmidt operator.
In summary,  both factors $J$ and $I_{k,\nu}$ in the  natural factorization \changeend
\begin{align*} 
\tridia {H_X}{\Lx 2 \nu}{H} {I_{X,\nu}}{J}{I_{k,\nu}} 
\end{align*} 
are Hilbert-Schmidt, and therefore their product $I_{X,\nu}$ is nuclear,
see e.g.~\cite[Chapter 11.3, Theorem 3]{BiSo87}.
\end{proofof}

\changebegin

\begin{proofof}{Corollary \ref{cor:general-result}}
Since $H\subset  \ell_\infty(T)$ we may consider the constant envelope function $b$ given by $b(t) = \snorm{\id:H\to  \ell_\infty(T) }$ for all $t\in T$.
Then all finite measures $\nu$ on $(T,\sB)$ satisfy $b\in \sLx 2 \nu$, and therefore the assertion follows from Theorem \ref{thm:general-result-envelope}.
\end{proofof}

\begin{proofof}{Corollary \ref{cor:restrict-to-bounded}}
As in the proof of Theorem \ref{thm:general-result-envelope} we may assume
without loss of generality 
that $H\subset \sLx 0 \sB$.

Now,
let $k$ be the kernel of $H$ and $k_{|T^*} := T^*\times T^* \to \R$
be its restriction onto $T^*$, that is, $k_{|T^*}(s,t) = k(s,t)$ for all $s,t\in T^*$.
It is well-known, see e.g.~\cite[Theorem 6 on page 25]{BeTA04}, that the RKHS of $k_{|T^*}$ is  
$H_{|T^*}:=\{ h_{|T^*}: h\in H  \}$ with norm given by 
\begin{align*}
\snorm f_{H_{|T^*}} = \inf\bigl\{ \snorm h_H: h\in H \mbox{ with } h_{|T^*} = f  \bigr\}
\end{align*}
for all $f\in H_{|T^*}$. Note that $H\subset \sLx 0 \sB$ immediately gives 
$H_{|T^*}\subset \sLx 0 {\sB_{|T^*}}$. 
In fact, even the inclusion $H_{|T^*}\subset \sLx \infty {\sB_{|T^*}}$ is true, since for 
$f\in H_{|T^*}$ and   $h\in H$ with 
$ h_{|T^*} = f$ we have 
\begin{align*}
\inorm{f} = \sup_{t\in T^*} |f(t)| = \sup_{t\in T^*} |h(t)| \leq \snorm h_H \cdot  \sup_{t\in T^*} b(t) \leq 
\snorm h_H \cdot M < \infty\, .
\end{align*} 
Since $\P(\{\sppath Y \in H\}) >0$ implies $H_X\subset H$ by Theorem \ref{thm:LuBe01a}, 
and this in turn guarantees $ H_{X_{|T^*}}\subset H_{|T^*}$, we then 
 also find the claimed inclusion
$ H_{X_{|T^*}}\subset \sLx \infty {\sB_{|T^*}}$.

 If we now write $Y^* :=  (Y_t)_{t\in T^*}$ for the restriction of $Y$ onto $T^*$, then 
$Y^*$ is a version of the centered Gaussian process $X_{|T^*}$.  Moreover,
restricting the paths of $Y$ onto $T^*$ gives the paths of $Y^*$, that is, 
\begin{align*}
\sppath Y^*(\om) = \bigl( \sppath Y(\om)\bigr)_{|T^*} \, , \myqquad \om\in \Om,
\end{align*}
and therefore we find $\P(\{ \sppath Y^* \in H_{|T^*}\}) \geq    \P(\{\sppath Y \in H\}) >0$. Applying  Corollary \ref{cor:general-result} to $X_{|T^*}$, $H_{|T^*}$, and the version $Y^* =(Y_t)_{t\in T^*}$ then gives the assertion.
\end{proofof}

\changeend


\subsection{Proofs Related to Essentially Equivalent RKHSs}

In this subsection we prove Lemma \ref{lem:kernel-diff} as well as Theorems
\ref{thm:ess-eq-rkhs} and \ref{thm:constr-essent-equiv-dom-rkhs}.
 
To begin with,  
let us assume that $H_1$ and $H_2$ are RKHSs on $T$ with kernels $k_1$ and $k_2$. Then their sum $H_1 + H_2$ 
can be equipped
with the natural norm 
\begin{align}\label{eq:sum-rkhs}
\snorm h_{H_1 + H_2}^2 = \inf\bigl\{ \snorm{h_1}_{H_1}^2 + \snorm {h_2}_{H_2}^2: h_1\in H_1, h_2\in H_2  \mbox{ with } h=h_1+h_2    \bigr\} \, .
\end{align}
It turns out that $H_1+H_2$ equipped with this norm is 
the RKHS of the kernel $k_1+k_2$, see e.g.~\cite[Theorem on page 353]{Aronszajn50a}.
Consequently, \eqref{eq:ess-equiv} implies $H_1 + \hstar 1 \cong H_2 + \hstar 2$ as noted at the beginning of Section \ref{sec:prelims}.
With these preparations, the proof of Lemma \ref{lem:kernel-diff} is straightforward:

\begin{proofof}{Lemma \ref{lem:kernel-diff}}
By our assumption we have $k_1 + k_1^\star = k_2 + k_2^\star$. As discussed above this gives $H_1 + \hstar 1 = H_2 + \hstar 2$, and hence 
also $H_1 + \hstar 1 \smspac H_2 + \hstar 2$.
\end{proofof}

In the following,  we adopt the usual notation  $H_1\oplus H_2 :\smspac H_1 + H_2$ 
if we additionally have $H_1 \cap H_2 = \{0\}$. Moreover, 
if we  equip $H_1\oplus H_2$ with the norm given by \eqref{eq:sum-rkhs}, then the resulting RKHS is denoted by 
$H_1\oplus_2 H_2$. 
In this case,
taking the infimum in \eqref{eq:sum-rkhs} is, of course, superfluous, and 
$H_1$ and $H_2$ are closed subspaces of $H_1\oplus_2 H_2$ with $H_1^\perp = H_2$ and $H_2^\perp = H_1$.

Now, our first auxiliary result shows that sums of RKHSs are  suitable direct sums if one of the involved spaces is 
finite dimensional.

\begin{lemma}\label{lem:essen-equiv-direct} 
Let $H$ and $H_0$ be RKHSs on $T$ with  $\dim H_0 < \infty$. Then there exists 
an RKHS $\hstar{}$ on $T$ with $\hstar{}\subset H_0$ and 
\begin{align}\label{lem:essen-equiv-direct-h0} 
H \oplus \hstar{} \smspac H+H_0  \, .
\end{align}
\end{lemma}

\begin{proofof}{Lemma \ref{lem:essen-equiv-direct}}
Since $H\cap H_0$ is a subspace of $H_0$ and $\dim H_0<\infty$, there exists
a subspace $\hstar{}\subset H_0$ with 
\begin{align}\label{lem:essen-equiv-direct-h1} 
H_0 \smspac (H\cap H_0) \oplus \hstar{}\, .
\end{align}
Moreover, $\dim H_0<\infty$ shows that $\hstar{}$ is closed in $H_0$ and thus $\hstar{}$ is an RKHS.

We first show that $H\cap \hstar{} = \{0\}$. This this end, we pick an $h\in H\cap \hstar{}$. By $\hstar{}\subset H_0$
we find $h\in H\cap H_0$ and hence we have $h\in (H\cap H_0) \cap \hstar{}$. By \eqref{lem:essen-equiv-direct-h1}  
we then obtain $h= 0$.

For the verification of \eqref{lem:essen-equiv-direct-h0}  we first note that $\hstar{}\subset H_0$ implies 
$H \oplus \hstar{} \subset H+H_0$.
For the converse inclusion, we fix some $h\in H$ and $h_0\in H_0$. By  \eqref{lem:essen-equiv-direct-h1} there then exist
a $g\in H\cap H_0$ and an   $h^\star\in \hstar{}$ with $h_0 = g + h^\star$. This gives 
\begin{align*}
h+h_0 = h + h_0-h^\star + h^\star       \in H + \hstar{} \smspac H \oplus \hstar{}\, ,
\end{align*}
where in the second step we used $h_0-h^\star = g \in H\cap H_0 \subset H$.
\end{proofof}

The next auxiliary  lemma shows that we may use direct sums in the definition of essentially equivalent RKHSs and that we have some 
additional control over the spaces $H_1^\star$ and $H_2^\star$.

\begin{lemma}\label{lem:essent-equiv-in-E}
Let $H_1$ and $H_2$ be essentially equivalent RKHSs on $T$. If $E$ \changebegin is a vector space consisting of functions $T\to \R$ \changeend  with $H_1\cup H_2 \subset E$, then there exist
finite dimensional 
RKHSs $\hstar 1$ and $\hstar 2$ on $T$ with $\hstar 1 \cup  \hstar 2 \subset E$ and 
\begin{align*} 
H_1 \oplus \hstar 1 \smspac H_2 \oplus \hstar 2\, .
\end{align*}
\end{lemma}

\begin{proofof}{Lemma \ref{lem:essent-equiv-in-E}}
By our assumption and Lemma \ref{lem:essen-equiv-direct} there exist finite dimensional RKHSs $\htim 1$ and $\htim 2$ on $T$ with 
\begin{align}\label{lem:essent-equiv-in-E-h1}
H_1 \oplus \htim 1 \smspac H_2 \oplus \htim 2\, .
\end{align}
We equip both direct sums with the $\oplus_2$-norm and consider the two orthogonal projections
$P_i :H_i\oplus_2 \htim i \to H_i\oplus_2 \htim i$ onto $H_i$ as well as the 
two orthogonal projections
$Q_i :H_i\oplus_2 \htim i \to H_i\oplus_2 \htim i$ onto $\htim i$.
Moreover, we define $H:\smspac H_1 \oplus \htim 1 \smspac H_2 \oplus \htim 2$ and the subspace
\begin{align*}
\hdag{}  :\smspac \bigl\{(Q_1 - Q_2) h: h\in H\bigr\} \smspac \bigl\{(Q_2 - Q_1) h: h\in H\bigr\} \, .
\end{align*}
Clearly, $\hdag{}$ is a subspace of $\htim 1 + \htim 2$ and thus a finite dimensional RKHS on $T$.
Moreover, for $h\in H$ 
our construction yields
\begin{align}\label{lem:essent-equiv-in-E-h2}
P_1 h + Q_1h = h = P_2 h + Q_2 h\, .
\end{align}
In particular, for $h_1\in H_1$ we find 
$h_1 = P_1 h_1 =  P_2 h_1 + (Q_2 - Q_1)h_1 \in H_2 + \hdag{} $, and hence $H_1\subset   H_2 + \hdag{}$.
The latter implies  $H_1 + \hdag{}   \subset H_2 + \hdag{}$. Since 
%
%
analogously we find $H_2 + \hdag{} \subset H_1 + \hdag{}$ by the symmetry of $\hdag{}$ in $Q_1$ and $Q_2$, we have 
\begin{align*}
H_1 + \hdag{} \smspac H_2 + \hdag{}\, .
\end{align*}
Moreover, \eqref{lem:essent-equiv-in-E-h2} also shows $(Q_1-Q_2)h = P_2 h - P_1 h \in H_2 + H_1 \subset E$, and hence we obtain
$\hdag {} \subset E$. Applying Lemma \ref{lem:essen-equiv-direct} on both sides then yields the assertions.
\end{proofof}

\begin{lemma}\label{lem:kolm-numb-ess-equiv}
Let $H_1, H_2$  and  $\hstar 1, \hstar 2$
be   RKHSs on $T$ with $H_1 \oplus \hstar 1 \smspac H_2 \oplus \hstar 2$ and $m_i := \dim \hstar i<\infty$.
For $i=1,2$ we write $\overline H_i := H_i \oplus_2 \hstar i$ and consider the  embedding maps $J_i:H_i \to \overline H_i$.
Finally, let $E$ be a Banach space and 
$S_i:\overline H_i \to E$ be bounded and linear operators with $S_1h=S_2 h$ for all $h\in \overline H_1 \smspac \overline H_2$.
Then there exists a constant $c\in (0,\infty)$ with 
\begin{align*}
d_{n+m_2}(S_1\circ J_1) &\leq c\,  d_n(S_2\circ J_2)\, , \\
d_{n+m_1}(S_2\circ J_2) &\leq c \, d_n(S_1\circ J_1)
\end{align*}
for all $n\geq 1$.
\end{lemma}

\begin{proofof}{Lemma \ref{lem:kolm-numb-ess-equiv}}
We first note that we have
$\overline H_1 \cong \overline H_2$ as mentioned at the beginning of Section \ref{sec:prelims}, and hence there exists a 
constant $c\in (0,\infty)$ with 
\begin{align}\label{lem:kolm-numb-ess-equiv-h1}
c^{-1}\, d_n(S_1)  \leq d_n(S_2) \leq c\, d_n(S_1)\, , \qquad \qquad n\geq 1.  
\end{align}
Let us further write $P_i:\overline H_i\to \overline H_i$ for the orthogonal projection 
onto $H_i$ and $P^\perp_i:\overline H_i\to \overline H_i$ for the orthogonal projection 
onto its complement $H_i^\perp = \hstar i$, that is $P^\perp_i = \id_{\overline H_i} - P_i$, 
where $\id_{\overline H_i}$ is the identity map on $\overline H_i$.
Moreover, 
we define $Q_i:\overline H_i\to H_i$ by $Q_i h := P_i h$ for all $h\in \overline H_i$.
In other words, the only difference between $P_i$ and $Q_i$ are the spaces the operators formally map into.
Now, by construction the operators  $Q_i$ are metric surjections and we have 
\begin{align*}
P_i = J_i \circ Q_i \qquad \qquad \mbox{ and } \qquad \qquad Q_i\circ J_i = \id_{H_i}\, .
\end{align*}
By the additivity of the Kolmogorov numbers  
we then find for all $n\geq 1$ 
\begin{align*}
d_{n+m_i}(S_i) 
\leq d_{m_i+1}(S_i-S_i\circ P_i) + d_n(S_i\circ P_i)
&= d_{m_i+1}(S_i \circ P^\perp_i) + d_n(S_i \circ J_i \circ Q_i) \\
&\leq \snorm {S_i}\cdot  d_{m_i+1}( P^\perp_i) + d_n(S_i \circ J_i \circ Q_i) \\
& = d_n(S_i \circ J_i)\, ,
\end{align*}
where in the last step we used $\rank  P_i^\perp = m_i<m_i+1$ and the surjectivity of the Kolmogorov numbers.
Combining this inequality with $\snorm{J_i}= 1$ and \eqref{lem:kolm-numb-ess-equiv-h1} yields 
\begin{align*}
d_{n+m_2}(S_1\circ J_1) \leq d_{n+m_2}(S_1) \leq c\,d_{n+m_2}(S_2) \leq  c\,d_n(S_2 \circ J_2)\, ,
\end{align*}
and the second inequality can be shown analogously. 
\end{proofof}

\begin{proofof}{Theorem \ref{thm:ess-eq-rkhs}}
Since $H_1$ and $H_2$  are essentially equivalent  with
 $H_1,H_2\subset \changebegin\sLx 2 \nu$ by Lemma \ref{lem:Iknu-new}\changeend, Lemma \ref{lem:essent-equiv-in-E} 
shows that there are finite dimensional
RKHSs $\hstar 1\subset \changebegin\sLx 2 \nu\changeend$ and $\hstar 2\subset\changebegin \sLx 2 \nu\changeend$  on $T$ 
with  
\begin{align*}
H_1 \oplus \hstar 1 \smspac H_2 \oplus \hstar 2\, .
\end{align*}
We write $m_i := \dim \hstar i$ and observe that 
$\overline H_i := H_i \oplus_2 \hstar i \subset \changebegin\sLx 2 \nu$ \changeend 
for $i=1,2$. Consequently, we can consider the 
``embedding'' operators $S_i := I_{\overline H_i,\nu} :\overline H_i \to \Lx 2 \nu$ \changebegin  given by $h\mapsto \eqclass h$. \changeend
Clearly this shows  $S_1h=S_2 h$ for all $h\in \overline H_1 \smspac \overline H_2$ and for 
 the  embedding maps $J_i:H_i \to \overline H_i$ we have 
 \begin{align*}
 S_i \circ J_i = I_{H_i,\nu}\, .
 \end{align*}
 \changebegin Before be can apply  Lemma \ref{lem:kolm-numb-ess-equiv} we 
 finally need to verify that both $S_1$ and $S_2$ are  continuous: To this end, we note that   both operators $I_{H_i,\nu} :H_i \to \Lx 2 \nu$ 
 are continuous by our assumptions and Lemma \ref{lem:Iknu-new}, and the operators 
 $I_{\hstar i,\nu} :\hstar i \to \Lx 2 \nu$ are continuous because of $\dim \hstar i<\infty$. 
 For $h \in \overline H_i := H_i \oplus_2 \hstar i$ and its unique decomposition $h=h_i + h_i^\star$ with $h_i\in H_i$ and $h_i^\star \in \hstar i$ we then find 
 \begin{align*}
 \snorm {S_i h  }_{\Lx 2 \nu} 
 \leq  \mnorm {\eqclass {h_i} }_{\Lx 2 \nu} + \mnorm {\eqclass{h_i^\star}}_{\Lx 2 \nu}
 &\leq \snorm{I_{H_i,\nu}} \cdot \snorm{h_i}_{H_i} +    \snorm{I_{\hstar i,\nu}} \cdot \snorm{h_i^\star}_{\hstar i} \\
 & \leq \sqrt 2 \cdot \max\{   \snorm{I_{H_i,\nu}},  \snorm{I_{\hstar i,\nu}}  \} \cdot \bigl(  \snorm{h_i}_{H_i}^2 +  \snorm{h_i^\star}_{\hstar i}^2 \bigr)^{1/2} \\
 & = \sqrt 2 \cdot \max\{   \snorm{I_{H_i,\nu}},  \snorm{I_{\hstar i,\nu}}  \} \cdot \snorm h_{\overline H_i}\, ,
 \end{align*}
 that is, the operators $S_i: \overline H_i \to \Lx 2 \nu$ are indeed continuous.
 Applying  \changeend Lemma \ref{lem:kolm-numb-ess-equiv} thus gives a constant $c\in (0,\infty)$ with 
\begin{align}\label{prop:ess-eq-rkhs-h1}
d_{n+m_2}(I_{H_1,\nu}) \leq c\,  d_n(I_{H_2,\nu}) 
\myqquad \mbox{ and } \myqquad 
d_{n+m_1}(I_{H_2,\nu}) \leq c \, d_n(I_{H_1,\nu})
\end{align}
for all $n\geq 1$.
Let us now assume that  $I_{H_1,\nu}$ is nuclear. Then  $I_{H_1,\nu}$ is compact and \eqref{eq:kol-vsr-sing} together with the discussion 
around it shows both 
$s_n(I_{H_1,\nu}) = d_n(I_{H_1,\nu})$ and $d_n(I_{H_1,\nu})\to 0$. 
By \eqref{prop:ess-eq-rkhs-h1}  we conclude that $d_n(I_{H_2,\nu})\to 0$ and hence $I_{H_2,\nu}$
is also compact as discussed around \eqref{eq:kol-vsr-sing}. Consequently,  we have 
$s_n(I_{H_2,\nu}) = d_n(I_{H_2,\nu})$, and therefore the second inequality in \eqref{prop:ess-eq-rkhs-h1}
shows 
\begin{align*}
s_{n+m_1}(I_{H_2,\nu}) \leq c \, s_n(I_{H_1,\nu})
\end{align*}
for all $n\geq 1$. Since we assumed that  $I_{H_1,\nu}$ is nuclear, we then see that  $I_{H_2,\nu}$ is also nuclear.
The converse implication can be shown analogously. 
\end{proofof}

\begin{proposition}\label{prop:nuk-doim-vs-ess-equiv}
Let $H_1, H_2$ be RKHS on $T$ and    $\hstar 1, \hstar 2$
be  finite dimensional RKHSs on $T$ with $H_1 \oplus \hstar 1 \smspac H_2 \oplus \hstar 2$. 
Moreover, let $H$ be an RKHS on $T$ with $H_1\ll H$. Then there exists an RKHS $\hdag 1\subset \hstar 1$ with 
\begin{align*}
H_2 \ll H\oplus_2 \hdag 1\, .
\end{align*}
\end{proposition}

\begin{proofof}{Proposition \ref{prop:nuk-doim-vs-ess-equiv}}
For $i=1,2$ we 
consider the  embedding maps $J_i:H_i \to H_i \oplus_2 \hstar i$.
Now, we first note that by Lemma \ref{lem:essen-equiv-direct}  there exists an $\hdag 1\subset \hstar 1$
with 
\begin{align*}
H\oplus \hdag 1 \smspac H+ \hstar 1\, .
\end{align*}
%
Moreover,  by our assumption  $H_1\ll H$    we have  $H_1\subset H$ and   the resulting 
embedding   $J:H_1 \to H$ is Hilbert-Schmidt. 
Since $H_1\subset H$ gives
\begin{align*}
H_2 \oplus \hstar 2 \smspac H_1 \oplus \hstar 1 \subset H + \hstar 1 \smspac H\oplus \hdag 1\, ,
\end{align*}
 we can thus consider the two resulting embedding maps $S_i: H_i \oplus_2 \hstar i \to H\oplus_2 \hdag 1$.
 Clearly, this gives $S_ih=S_2 h$ for all $h\in H_1 \oplus_2 \hstar 1 \smspac H_2 \oplus_2 \hstar 2$. Moreover, the operators $S_i$ are continuous 
 since embedding maps between RKHSs are automatically continuous by a 
 simple application of the closed graph theorem, see e.g.~\cite[Lemma 23]{ScSt25a}. 
 By Lemma \ref{lem:kolm-numb-ess-equiv} there then exists a constant $c\in (0,\infty)$ such that for $m_1:= \dim H_1^\dagger$  we have 
\begin{align*}
d_{n+m_1}(S_2\circ J_2) &\leq c \, d_n(S_1\circ J_1)\, , \myqquad n\geq 1.
\end{align*}
Let us now consider the embedding map $I:H\to H\oplus_2 \hdag 1$, which by construction 
satisfies $\snorm I\leq 1$.
Clearly, our construction gives 
\begin{align*}
\quadiass {H_1}{H_1 \oplus_2 \hstar 1}H {H\oplus_2 \hdag 1} {J_1}J{S_1}I 
\end{align*}
and hence we obtain 
\begin{align*}
d_{n+m_1}(S_2\circ J_2) 
\leq c \, d_n(S_1\circ J_1)
= c \, d_n(I\circ J)
\leq  c  \, d_n(J)
\end{align*}
for all $n\geq 1$. Now, $J$ is Hilbert-Schmidt, and hence compact. This shows both 
$d_n(J) \to 0$ and 
$d_n(J) = s_n(J)$  
 as 
discussed around \eqref{eq:kol-vsr-sing}.
Consequently, we find $d_{n}(S_2\circ J_2)\to 0$ and thus also $s_{n}(S_2\circ J_2) = d_{n}(S_2\circ J_2)$, where for the second 
assertion we again used 
 the properties mentioned around  \eqref{eq:kol-vsr-sing}.
In summary, we thus have  $s_{n+m_1}(S_2\circ J_2) \leq  c  \, s_n(J)$, and since $J$ is Hilbert-Schmidt, so ist 
 $S_2\circ J_2$. However, $S_2\circ J_2 : H_2 \to H\oplus_2 \hdag 1$ is the embedding map, and hence we have found  $H_2 \ll H\oplus_2 \hdag 1$.
\end{proofof}

\begin{proofof}{Theorem \ref{thm:constr-essent-equiv-dom-rkhs}}
Let us consider the RKHS $E:= H^\dagger+ H_X$. 
By our assumptions we have $H_X\cup H \subset H_X\cup H^\dagger\subset E$, and therefore Lemma \ref{lem:essent-equiv-in-E}
provides us with finite dimensional 
 RKHSs $\hstar {}$ and $\hstar X$ on $T$ with $\hstar {} \cup \hstar X\subset E$ and 
\begin{align*} 
H \oplus \hstar {} \smspac H_X \oplus \hstar X\, .
\end{align*}
Combining this with the assumed $H\ll H^\dagger$, Proposition \ref{prop:nuk-doim-vs-ess-equiv} gives us an
 RKHS $H^\ddagger\subset \hstar {}$ with 
\begin{align*}
H_X \ll H^\dagger \oplus_2 H^\ddagger\, .
\end{align*}
In addition, note that $H_X^\dagger:=  H^\dagger \oplus_2 H^\ddagger$ is essentially equivalent to $H^\dagger$
and our construction ensures 
$H_X^\dagger=  H^\dagger \oplus_2 H^\ddagger\subset H^\dagger +   \hstar {} \subset H^\dagger +  H^\dagger+ H_X =   H^\dagger+ H_X$,
where for the second inclusion we used $H^*\subset E= H^\dagger+ H_X$. 
\end{proofof}


\subsection{Proof of Theorem \ref{thm:ew-char}}

\begin{proofof}{Theorem \ref{thm:ew-char}}
\ada i \changebegin We first recall that
Lemma \ref{lem:Iknu-new}
ensures the compactness of $I_{X,\nu}$ as 
 since $b\in \sLx 2 \nu$ is an envelope function for $H_X$. \changeend
%
The definition of singular numbers then yields
\begin{align*}
s_i( I_{X,\nu}) = \sqrt{\mu_i(I_{X,\nu} \circ  I_{X,\nu}^*)  } =  \sqrt{\mu_i(T_{X,\nu})  }\, .
\end{align*}
By \eqref{eq:nuc-op} we thus obtain the equivalence.

\ada {ii} Let $J:H_X\to H_X^\b$ be the embedding map. For $i\in I$ we then have 
\begin{align*}
\snorm{J (\sqrt{\mu_i} e_i)}_{H_X^\b} = \mu_i^{1/2 - \b/2}  \snorm{ \mu_i^{\b/2} e_i}_{H_X^\b} = \mu_i^{1/2 - \b/2} \, ,
\end{align*}
and this in turn shows 
\begin{align*}
\sum_{i\in I} \snorm{J (\sqrt{\mu_i} e_i)}_{H_X^\b}^2 = \sum_{i\in I}  \mu_i^{1-\b} \, .
\end{align*}
Since the $H_X$-positivity of $\nu$ ensures that 
$(\sqrt{\mu_i} e_i)_{i\in I}$ is an ONB of $H_X$, see \cite[Theorem 3.1]{StSc12a}, we can then conclude 
that \eqref{eq:nuc-int-op-dominance} indeed characterizes the Hilbert-Schmidt property of $J$, see the 
discussion around \eqref{eq:hs-op}.
The second assertion directly follows from $H_X\ll H_X^{\b}$ and Theorem \ref{thm:LuBe01a}.  
%
%
%
%
\end{proofof}


\subsection{Proofs Related to RKHSs of  Sobolev Type}\label{sec:proofs-hoel+sobol}

Before we can prove Theorem \ref{thm:sobol-main} some preparations are required. We begin by providing details for 
\eqref{eq:bessel-rkhs-short}. To this end, let 
 $B^s_{p,q}(\Rd)$ and $F^s_{p,q}(\Rd)$   denote the  Besov spaces, respectively  Triebel-Lizorkin spaces, see e.g.~\cite{AdFo03,EdTr96}.
 Then  
 for $s>d/2$  we have 
\begin{align}\label{eq:bessel-rkhs}
 H^s(\R^d) \cong F_{2,2}^s(\Rd) \cong  B_{2,2}^s(\Rd) \hookrightarrow  C_b(\Rd)\, ,
\end{align}
where the first identity follows from  e.g.~\cite[Chapter 7, Paragraphs 7.62, 7.65, and 7.57]{AdFo03},
the second identity is shown in e.g.~\cite[Chapter 7, Paragraph 7.67]{AdFo03}, and the inclusion can be found in 
 e.g.~\cite[Theorem 7.34]{AdFo03}. 
 
In addition, Triebel-Lizorkin and Besov spaces on open Euclidean balls $T\subset \Rd$ are defined analogously to \eqref{def:HsT}, that is 
$A^s_{p,q}(T) := \{f_{|T} :    f\in A^s_{p,q}(\Rd)\}$ with norm
\begin{align*}
\snorm{h}_{A^s_{p,q}(T)} := \inf\bigl\{ \snorm f_{A^s_{p,q}(\Rd)}: f\in A^s_{p,q}(\Rd) \mbox{ with } f_{|T} = h  \bigr\}\,
\end{align*}
where $A^s_{p,q}$ stands for either $F^s_{p,q}$ or $B^s_{p,q}$, see e.g.~\cite[Chapter 2.5]{EdTr96}.
Note that if $F_{p,q}^s(\Rd)$ or  $B_{p,q}^s(\Rd)$
are   BSFs then the same is true for $F_{p,q}^s(T)$ or  $B_{p,q}^s(T)$. In addition, continuous inclusions between such spaces on $\Rd$
carry over to the restricted spaces on $T$, see e.g.~\cite[Chapter 2.5]{EdTr96} or \cite[Proposition 3]{ScSt25a} for a general BSF statement.
For $s>d/2$,  Equation \eqref{eq:bessel-rkhs} thus implies 
\begin{align}\label{eq:bessel-rkhs-onT}
 H^s(T) \cong F_{2,2}^s(T) \cong  B_{2,2}^s(T) \hookrightarrow  C_b(T)\, .
\end{align} 


For the proof of Theorem \ref{thm:sobol-main} we need the following auxiliary result that
investigates both the operators $I_{H^s(T),\lbd}: H^s(T)\to \Lx 2 T$
and the nuclear dominance between fractional Sobolev spaces.

\begin{lemma}\label{nuclear-stuff-for-sobol}
 Let $T\subset \Rd$ be bounded and measurable with non-empty interior and $s>d/2$. Then the following statements hold true:
 \begin{enumerate}
 \item The operator $I_{H^s(T),\lbd}: H^s(T)\to \Lx 2 T$ is nuclear if and only if $s>d$. 
 \item If $s>d$ and $r>d/2$, then $H^s(T) \ll H^r(T)$ if and only if $r<s-d/2$.
 \end{enumerate}
\end{lemma}

\begin{proofof}{Lemma \ref{nuclear-stuff-for-sobol}}
\ada i
 We first consider the case, when $T$ is an open Euclidean ball. Then we have $B_{2,2}^0(\Rd) \cong F_{2,2}^0(\Rd)\cong \Lx 2 \Rd$, 
 where the first identity follows from 
  e.g.~\cite[(1) in Chapter 2.3.3]{EdTr96}  and the second identity can be found in e.g.~\cite[(5) in Chapter 2.5.6]{Triebel83}. 
  Restricting to $T$ leads to  $B_{2,2}^0(T) \cong \Lx 2 T$, and therefore 
   \cite[(10) in Chapter 2.5.1]{EdTr96} 
 shows that the embedding $B_{2,2}^s(T) \hookrightarrow \Lx 2 T$
 is compact.
 By \eqref{eq:bessel-rkhs-onT}  we can thus  conclude that  $I_{H^s(T),\lbd}: H^s(T)\to \Lx 2 T$ is also compact, and by 
%
%
%
%
\eqref{eq:kol-vsr-sing} we hence find
\begin{align}\label{ex:ou-1-h1}
a_n(I_{H^s(T),\lbd}) = d_n(I_{H^s(T),\lbd}) =  s_n(I_{H^s(T),\lbd})\, , \myqquad n\geq 1.
\end{align}
Moreover,  the already observed $H^s(T) \cong B_{2,2}^s(T)$ and $B_{2,2}^0(T) \cong   \Lx 2 T$ in combination with the results 
%
%
in 
\cite[Chapter 3.3.4]{EdTr96} show that 
\begin{align}\label{ex:ou-1-h2}
a_n(I_{H^s(T),\lbd})  \sim n^{-s/d} \, .
\end{align}
Combining \eqref{ex:ou-1-h1} with \eqref{ex:ou-1-h2}  then yields the characterization.

To prepare the proof of  the general case, let us assume that we have measurable $T_1\subset T_2 \subset \Rd$ with 
$\lbd (T_1), \lbd(T_2) \in (0,\infty)$. Moreover, let $Q:H^s(T_2) \to H^s(T_1)$ and $R:\Lx 2 {T_2} \to \Lx 2 {T_1}$ be 
the restriction operators. Then $Q$ and $R$ are metric surjections and we have the commutative diagram
\begin{align*}
\quadiass {H^s(T_2)}{\Lx 2 {T_2}}{H^s(T_1)} {\Lx 2 {T_1}} {I_{H^s(T_2),\lbd}} QR {I_{H^s(T_1),\lbd}} \, .
\end{align*}
By the surjectivity of the Kolmogorov numbers we thus find 
\begin{align}\label{ex:ou-1-hxx}
d_n(I_{H^s(T_1),\lbd}) 
= 
d_n(I_{H^s(T_1),\lbd} \circ Q) 
= 
d_n (R \circ I_{H^s(T_2),\lbd}) 
\leq 
d_n(I_{H^s(T_2),\lbd})\, , \myqquad n\geq 1.
\end{align}

With these preparations we can now establish the assertion in the general case. To this end, we first note that 
by our assumptions on $T$ 
there exist open balls $B_0\subset T \subset B_1$. Now, we have already seen at the very
beginning that $I_{ H^s(B_i),\lbd}$ is compact, and hence we have $d_n(I_{H^s(B_i),\lbd}) \to 0$ for $i=0,1$.
Applying \eqref{ex:ou-1-hxx} we further find $d_n(I_{H^s(T),\lbd}) \to 0$, and therefore 
$I_{H^s(T),\lbd}$ is also compact. By \eqref{eq:kol-vsr-sing} we conclude that 
\begin{align*}
s_n(I_{H^s(T),\lbd}) = d_n(I_{H^s(T),\lbd})\, , \myqquad n\geq 1,
\end{align*}
while our initial considerations \eqref{ex:ou-1-h1} with \eqref{ex:ou-1-h2}  showed 
\begin{align*}
d_n(I_{H^s(B_i),\lbd}) = a_n(I_{H^s(B_i),\lbd})  \sim n^{-s/d}\, , \myqquad n\geq 1
\end{align*}
for $i=0,1$. Combining these identities with \eqref{ex:ou-1-hxx} yields $s_n(I_{H^s(T),\lbd}) \sim  n^{-s/d}$.

\ada {ii} Again, we first consider the case, when $T$ is an open ball. Here, we note that  $H^s(T) \cong B_{2,2}^s(T)$ 
and  $H^r(T) \cong B_{2,2}^r(T)$ in combination with  
the results in 
\cite[Chapter 3.3.4]{EdTr96} show that 
\begin{align*} 
a_n(J )  \sim n^{-(s-r)/d}\, ,
\end{align*}
where $J: H^s(T) \to H^r(T)$ denotes the embedding map. By the results mentioned around \eqref{eq:kol-vsr-sing} the same is 
true for the Kolmogorov numbers and in particular we have $d_n(J) \to 0$. Therefore, $J$ is compact and 
 \eqref{eq:kol-vsr-sing} shows that 
 \begin{align*}
 s_n(J )  \sim n^{-(s-r)/d}\, .
 \end{align*}
 Consequently, we have $\sum_{n\geq 1} s_n^2(J) < \infty$ if and only if $r< s-d/2$.
%

To prepare the proof of  the general case, we again assume that we have measurable $T_1\subset T_2 \subset \Rd$ with 
$\lbd (T_1), \lbd(T_2) \in (0,\infty)$. 
Moreover, let $Q_s:H^s(T_2) \to H^s(T_1)$ and $Q_r:H^r(T_2) \to H^r(T_1)$ be 
the restriction operators, which by construction of the involved spaces are metric surjections.
In addition, note that $H^s(\Rd) \subset H^r(\Rd)$ ensures that the embeddings $J_i: H^s(T_i) \to H^r(T_i)$
are defined for $i=1,2$. Clearly, these maps enjoy the following commutative diagram
\begin{align*}
\quadiass {H^s(T_2)}{H^r(T_2)}{H^s(T_1)} {H^r(T_1)} {J_2} {Q_s}{Q_r} {J_1} \, .
\end{align*}
By the surjectivity of the Kolmogorov numbers we thus find 
\begin{align*}
d_n(J_1) 
= 
d_n(J_1 \circ Q_s) 
= 
d_n (Q_r \circ J_2) 
\leq 
d_n(J_2)\, , \myqquad n\geq 1.
\end{align*}
With these preparations, the rest of the proof follows the arguments used in the proof of the general case of \emph{i)}.
\end{proofof}

\changebegin

In addition, the proof of Theorem \ref{thm:sobol-main}
as well as some of the examples in Section \ref{sec:examples} 
need the following result that 
refines the analysis of Theorem \ref{thm:nice-rkhs}.

\begin{lemma}\label{lem:nice-bounded-int}
Let $T\subset \Rd$ be bounded and measurable with non-empty interior. Moreover, let 
$(\Om,\sA,\P)$ be a complete probability space and $X:= (X_t)_{t\in T}$ be a centered Gaussian process over  $(\Om,\sA,\P)$
whose  RKHS $H_X$  consists of continuous functions. Finally, let 
$H$ be an RKHS  $H$ on $T$  
and $Y$ be a version   of $X$ with 
\begin{align*}
\P(\{\sppath Y \in H\}) > 0\, .
\end{align*}
Then there  exist a version $Y^*$ of $X$, a  separable RKHS $H^*$ on $T$ with kernel $k^*$, and an $M>0$ such that $H_X\subset H^*$ is dense in $H^*$,
the set $T^* := \{t\in T: k(t,t) \leq M\}$ is measurable and has non-empty interior, and 
\begin{align*}
\P(\{\sppath Y^* \in H^*\}) =1 \, .
\end{align*}
\end{lemma}

\begin{proofof}{Lemma \ref{lem:nice-bounded-int}}
By Theorem \ref{thm:nice-rkhs}  we may assume without loss 
of generality that $H$ is separable and 
$H_X\subset H$ is dense in $H$. In the following we thus show the assertions
for $H$ and $Y$.

Under these assumptions, the embedding $J:H_X\to H$   satisfies
\begin{align*}
\overline {\ran J}^{H} = \overline {H_X}^H = H\, , 
\end{align*}
and therefore Lemma \ref{lem:comput-dominance} gives us an at most countable 
family $(e_i) \subset H_X$ such that the functions $f_i := Je_i = e_i$ 
form an ONB $(f_i)_{i\in I}$ of $H$. By 
e.g.~\cite[Theorem 4.20]{StCh08} we can thus conclude  that the kernel $k$ of $H$
is given by 
\begin{align*}
k(s,t) = \sum_{i\in I} f_i(s) f_i(t)\, , \myqquad s,t\in T.
\end{align*}
Now, since $H_X$ consists of continuous functions, all $f_i$ are continuous, and 
therefore the function $b:T\to [0,\infty)$ defined by $b(t) :=  {k(t,t)}$ is both measurable and
lower semi-continuous. 

We now fix an  open and non-empty $O\subset T$ and 
write
   $O_M := \{t\in O: b_{|O}(t) > M\}$. Obviously, this yields
 \begin{align}\label{eq:baire-ansatz}
\emptyset = \bigcap_{M\in \N} O_M 
 \end{align} 
 and all sets $O_M$ are relatively open in $O$ by the  lower 
 semi-continuity of $b$. Moreover, 
 \eqref{eq:baire-ansatz} shows that their countable intersection is not dense in $O$.
In addition, since $\Rd$ is a complete metric space, 
 there exists a \emph{complete} metric on $O$ describing the topology  
 on $O$, see e.g.~\cite[Theorem 2.5.4]{Dudley02}.
Consequently, Baire's theorem, see e.g.~\cite[Theorem 1.5.4]{Megginson98}, shows that 
there is an $M\geq 1$ such that $O_M$ is not dense in $O$.
This in turn means that its complement $O\setminus O_M = \{t\in O: b_{|O}(t) \leq M\}$  
has non-empty relative  interior   in $O$. Consequently, there exists an relatively $O$-open set $O^*\neq \emptyset$ with 
\begin{align*}
O^* \subset \{t\in O: b_{|O}(t) \leq M\} \subset \{t\in T: b(t) \leq M\}\,.
\end{align*}
Since $O$ open in $\Rd$, so is $O^*$. In other words, the set $T^*$ has non-empty interior.
\end{proofof}

\changeend

\begin{proofof}{Theorem \ref{thm:sobol-main}}
\changebegin \atob {ii} {iii} Trivial. 

 \atob {iii} i By Lemma \ref{lem:nice-bounded-int} we may assume without loss 
 of generality that there is an $M>0$ such that 
 $T^* := \{t\in T: k(t,t) \leq M\}$ is measurable and has non-empty interior,
 where $k$ is the kernel of $H$.
 Consequently, 
 the Lebesgue measure on $T^*$ is non-trivial and finite.
By   Corollary \ref{cor:restrict-to-bounded}
we then  know that 
\begin{align*}
I_{X_{|T^*},\lbd}:H_{X_{|T^*}}\to \Lx 2 {T^*}
\end{align*}
 is nuclear. Moreover, since  $H_X$ is essentially equivalent to $H^s(T)$, the spaces 
 $H_{X_{|T^*}}$ and $H^s(T^*)$ are also essentially equivalent. 
 In addition, we know  $H_{X_{|T^*}} \subset \sLx \infty {T^*}$ by Corollary
 \ref{cor:restrict-to-bounded}, and since we also have 
 $H^s(T^*) \subset \sLx \infty {T^*}$ both spaces have a joint 
 constant envelope.
Consequently,
Theorem \ref{thm:ess-eq-rkhs}   shows that $I_{H^s(T^*),\lbd}: H^s(T^*)\to \Lx 2 {T^*}$ is also nuclear. Now the assertion
follows by Lemma \ref{nuclear-stuff-for-sobol}.

 \changeend

\atob i {ii} Let us fix an $d/2<r<s-d/2$. Lemma \ref{nuclear-stuff-for-sobol} then shows $H^s(T) \ll H^r(T)$.
Since $H_X$ is essentially equivalent to $H^s(T)$ by assumption, Theorem \ref{thm:constr-essent-equiv-dom-rkhs}
then provides an RKHS $H_X^\dagger$ on $T$ that is essentially equivalent to $H^r(T)$ and satisfies $H_X \ll H_X^\dagger$.
Now the assertion follows from Theorem \ref{thm:LuBe01a} for $H:= H_X^\dagger$.
\end{proofof}

\begin{proofof}{Theorem \ref{thm:mixed-sobol}}
Using \eqref{eq:bessel-rkhs-onT} we first note that it suffices to consider the case $H_0 = B^s_{2,2}(0,1)$. In the following we thus assume that 
 $\hat k_0$ is the   kernel  of
 $\hat H_0 := B_{2,2}^s(\R)$ and  that $k_0$ is its restriction to $(0,1)$.
 
Now let $S_{2,2}^sB(\Rd) $ denote the Besov space of dominating mixed smoothness $s$, see \cite[Appendix A.2]{SiUl09a} and the references mentioned therein.
Then \cite[Theorem 2.2]{SiUl09a}   shows 
\begin{align}\label{thm:mixed-sobol-h1}
S_{2,2}^sB(\Rd) \cong B^s_{2,2}(\R) \otimes_{\a_2}\dots  \otimes_{\a_2} B^s_{2,2}(\R)\, ,
\end{align}
where $\a_2$ denotes the $2$-nuclear tensor norm, see  \cite[Definition B.2]{SiUl09a}. Moreover, the $\a_2$-tensor norm equals the 
$g_2 :=\a_{2,1}$-tensor norm in \cite[Chapters 12.5 and 12.7]{DeFl93}, and $g_2$ in turn equals the natural Hilbert space tensor norm, see 
\cite[Chapters 26.6 and 26.7]{DeFl93}, where the latter is called Hilbert-Schmidt norm. Denoting the  natural Hilbert space tensor product by $\otimes_2$, 
the Identity \eqref{thm:mixed-sobol-h1} thus reads as 
\begin{align}\label{thm:mixed-sobol-h2}
S_{2,2}^sB(\Rd) \cong B^s_{2,2}(\R) \otimes_{2}\dots  \otimes_{2} B^s_{2,2}(\R)\, .
\end{align}
Now the tensor product space on the right hand side can be interpreted as  the RKHS of the product kernel 
\begin{align*}
\prod_{i=1}^d \hat k_0(t_i,t_i')\, , \myqquad t,t'\in \Rd,
\end{align*}
see e.g.~\cite[Chapter 1.4.6]{BeTA04} or \cite[Lemma 4.6]{StCh08}. Restricting both sides of \eqref{thm:mixed-sobol-h2} to $T$ thus gives 
\begin{align}\label{thm:mixed-sobol-h3}
S_{2,2}^sB(T) \cong B^s_{2,2}(0,1) \otimes_{2}\dots  \otimes_{2} B^s_{2,2}(0,1) = H_X\, ,
\end{align}
where the space $S_{2,2}^sB(T)$  is defined by restricting the functions of $S_{2,2}^sB(\Rd)$ to $T$ as in  
\cite[Definition 3.8 and Remark 3.9]{Nguyen17a}.
Moreover,   the embedding $I:S_{2,2}^sB(T) \to \Lx 2 T$ is defined and compact and its approximation numbers satisfy 
\begin{align}\label{thm:mixed-sobol-h4}
a_n(I) \sim n^{-s} (\ln n)^{(d-1)s}
\end{align}
 see \cite{Romanyuk01a} as well as \cite[Theorem 2.11]{Nguyen17a} for a presentation of the results 
  of \cite{Romanyuk01a} 
 in the notation
 we use in this proof.
 
 \atob {ii}i By \changebegin Corollary \ref{cor:general-result} \changeend we know that $I_{X,\lbd}:H_X\to \Lx 2 T$ is nuclear. 
 By combining \eqref{thm:mixed-sobol-h3}, \eqref{thm:mixed-sobol-h4}, and \eqref{eq:kol-vsr-sing}
we thus find $s>1$.
 
\atob i {ii} We pick an $r\in (1/2,s-1/2)$ and an ONB $(e_i)$ of $B^s_{2,2}(0,1)$. By Lemma \ref{nuclear-stuff-for-sobol} and \eqref{eq:bessel-rkhs-onT} we then know that 
$B^s_{2,2}(0,1) \ll B^r_{2,2}(0,1)$, and hence we have 
\begin{align*}
\sum_{i=1}^\infty \snorm{e_i}_{B^r_{2,2}(0,1)}^2 < \infty
\end{align*}
as discussed around \eqref{eq:hs-op}.
 Now recall that the collection of elementary tensors 
 $e_{i_1}\otimes \dots\otimes e_{i_d}$ forms an ONB of $B^s_{2,2}(0,1) \otimes_{2}\dots  \otimes_{2} B^s_{2,2}(0,1) = H_X$.
 Moreover, we have 
 \begin{align*}
 \sum_{i_1=1}^\infty \dots \sum_{i_d=1}^\infty \snorm{e_{i_1}\otimes \dots\otimes e_{i_d}}_{B^r_{2,2}(0,1) \otimes_{2}\dots  \otimes_{2} B^r_{2,2}(0,1)}^2
 &= 
 \sum_{i_1=1}^\infty \dots \sum_{i_d=1}^\infty \prod_{j=1}^d \snorm{e_{i_j}}_{B^r_{2,2}(0,1)}^2 \\
 &= 
 \biggl(  \sum_{i=1}^\infty \snorm{e_i}_{B^r_{2,2}(0,1)}^2 \biggr)^d \\
 &< \infty\, , 
 \end{align*}
 which shows that the embedding $B^s_{2,2}(0,1) \otimes_{2}\dots  \otimes_{2} B^s_{2,2}(0,1)\hookrightarrow B^r_{2,2}(0,1) \otimes_{2}\dots  \otimes_{2} B^r_{2,2}(0,1)$ is Hilbert-Schmidt. Now \emph{ii)} follows by Theorem \ref{thm:LuBe01a}.
\end{proofof}


\subsection{Proofs for the Examples}\label{sec:proofs-ex}

\begin{proofof}{Example \ref{ex:wp}} \changebegin 
It is well-known, see e.g.~\cite[Chapter III.2, page 69]{Adler90}, \cite[Example 1, page 71]{BeTA04}, or 
\cite[page 28]{Lifshits12},  that the RKHS $H_W$ of the Wiener process is given 
\begin{align}\label{eq:rkhs-wp}
H_W \smspac \{h\in H^1([0,1]): h=0\}\, .
\end{align}
Therefore, $H_W$ is essentially equivalent to $H^1([0,1])$ and the non-existence of $H$  follows from Theorem \ref{thm:sobol-main}.

For the proof of the second assertion  we note that $I_{H^1([0,1],\lb}: H^1([0,1]\to \Lx 2 \lb$ is not nuclear by Lemma \ref{nuclear-stuff-for-sobol}. Moreover, we obviously have $H_W \subset H^1([0,1])\subset \sLx \infty T$, 
and therefore Theorem \ref{thm:ess-eq-rkhs}  shows 
that $I_{W,\lb}: H_W\to \Lx 2 \lb$ is not nuclear, either. \changeend
%
\end{proofof}

\begin{proofof}{Example \ref{ex:bb}}
Let us recall that the covariance function of the Brownian bridge $X:= (X_t)_{t\in T}$ is 
\begin{align*}
k_{ X}(t_1,t_2) =  \min\{t_1,t_2\}-t_1t_2 =    k_W(t_1,t_2) - t_1t_2    \, , \myqquad t_1,t_2\in [0,1],
\end{align*}
where $k_W$ is the covariance function of the Wiener process.
Now,  $(t_1,t_2)\mapsto t_1t_2$ is a kernel on $[0,1]$ whose RKHS is one dimensional, and 
consequently, $H_B$ and $H_W$ are essentially equivalent by Lemma \ref{lem:kernel-diff}. \changebegin Moreover, we know by \eqref{eq:rkhs-wp} that 
$H_W$ is essentially equivalent to  $H^1([0,1])$, and therefore $H_B$ is also 
essentially equivalent to   $H^1([0,1])$. Now the assertion
 follows from Theorem \ref{thm:sobol-main} since $H_B$ consists of continuous functions. \changeend
%
%
%
\end{proofof}

\begin{proofof}{Example \ref{ex:ou-1}}
We first consider the covariance function $k_X^{(1)}$. For the corresponding RKHS $H_X^{(1)}$ we then find 
$H_X^{(1)} \cong H^1(T)$ by  e.g.~\cite[Example 17 on p.~316]{BeTA04} applied  for some $a,b\in \R$ with $T\subset (a,b)$ followed by a restriction onto $T$.
By Theorem \ref{thm:sobol-main}
we then  obtain the assertion.

In the case of the second covariance, we first observe that $(t_1,t_2)\mapsto a\exp(-\s(t_1+t_2))$ gives  a kernel $k^*$ on $[0,1]$ whose RKHS is one dimensional. Moreover, we obviously have $k_X^{(2)} = k_X^{(1)} - k^*$, and therefore  
the RKHS $H_X^{(2)}$ of $k_X^{(2)}$ is essentially equivalent to $H_X^{(1)}\cong  H^1(T)$ by
 Lemma \ref{lem:kernel-diff}. Now we again obtain the assertion by Theorem \ref{thm:sobol-main}.
\end{proofof}


\begin{proofof}{Example \ref{ex:fbm}}
In the following we write 
$H_\a := H_{B^{(\a)}}$ for the RKHS of $B^{(\a)}$. Then \cite[Theorem 6.12]{Picard11a} shows that a function $h:T\to \R$ is in 
$H_\a$ if and only if there are sequences $(a_n)_{n\geq 0}$ and $(b_n)_{n\geq 1}$ with 
\begin{align*}
\sum_{n\geq 0} n^{2\a+1} a_n^2 < \infty \qquad \mbox{ and } \qquad \sum_{n\geq 1} n^{2\a+1} b_n^2 < \infty
\end{align*}
such that for all $t\in [0,1]$ we have 
\begin{align*}
h(t) = h(1) t + \sum_{n\geq 0} a_n \cos(2\pi n t) + \sum_{n\geq 1} a_n \sin(2\pi n t)  \, .
\end{align*}
With these preparations we can now turn to the actual proof. 

Here, we first note that in the case $\a=1/2$ we have $B^{(\a)} = W$, and hence
there is nothing left to prove by Example \ref{ex:wp}.

Moreover, in the case $0<\a<1/2$ the above description of $H_\a$ yields $H_{1/2} \subset H_\a$ and the corresponding embedding 
map $J:H_{1/2}\to H_\a$ is continuous by a simple application of the closed graph theorem, see   \cite[Lemma 23]{ScSt25a}.
\changebegin
Now assume that there was an  RKHS $H$ on $T$  with 
$\P(\{\sppath B^{(\a)} \in H\}) =1$. By 
Lemma \ref{lem:nice-bounded-int} we may then assume without loss 
 of generality that there is an $M>0$ such that 
 $T^* := \{t\in T: k(t,t) \leq M\}$ is measurable and has non-empty interval,
 where $k$ is the kernel of $H$.
 Corollary \ref{cor:restrict-to-bounded} then shows
that the operator $I_{(H_\a)_{|T^*},\lb}:(H_\a)_{|T^*}\to \Lx 2 {T^*}$ is nuclear.
\changeend
%
%
The natural factorization 
\begin{align*} 
\tridia {(H_{1/2})_{\changed{|T^*} }}{\Lx 2 {T\changed{^*}}}{(H_\a)_{\changed{|T^*} }} {I_{(H_{1/2})_{\changed{|T^*} },\lb}}{J}{I_{(H_\a)_{\changed{|T^*} },\lb}} 
\end{align*} 
then shows that $I_{(H_{1/2})_{\changed{|T^*} },\lb}$ would also be nuclear. 
\changebegin 
However, we have seen in \eqref{eq:rkhs-wp} that $H_{1/2} = H_W$ is 
essentially equivalent to $H^1([0,1])$ and therefore 
$(H_{1/2})_{\changed{|T^*}}$ is essentially equivalent to $H^1(T^*)$. 
Now, Lemma \ref{nuclear-stuff-for-sobol} shows that 
$I_{H^1(T^*),\lb}: H^1(T^*)\to \Lx 2 {T^*}$ is not nuclear, and 
since we have both $H^1(T^*) \subset \sLx \infty {T^*}$ 
and $(H_{1/2})_{\changed{|T^*}}\subset \sLx \infty {T^*}$, 
Theorem \ref{thm:ess-eq-rkhs} shows that $I_{(H_{1/2})_{\changed{|T^*} },\lb}$ 
is not nuclear, either. In other words, we have found a contradiction.
\changeend

Finally, in the case $1/2<\a<1$ it is well known that $\P$-almost all paths of $B^{(\a)}$ are $\b$-H\"older continuous for all $\b\in (0,\a)$, where we note that 
this can be quickly verified as the stationary increments and the self-similarity of $B^{(\a)}$ gives
\begin{align*}
 \E \bigl| B^{(\a)}_{t_1} - B^{(\a)}_{t_2}  \big|^p =  \E \bigl| B^{(\a)}_{|t_1-t_2|}   \big|^p = |t_1-t_2|^{\a p } \cdot\E  \bigl| B^{(\a)}_{1}   \big|^p \, , \myqquad t_1,t_2\geq 0, p>0\, ,
\end{align*}
and therefore an application of the Kolmogorov-Chentsov theorem, see e.g.~\cite[Theorem 21.6]{Klenke14} for $p >0$  with $\b< (\a p-1)/p$ gives the $\b$-H\"older continuity.
Let us now fix some $s\in (1/2,\a)$ and some $\b\in (s,\a)$.
 By  e.g.~\cite[(6) and (9) in Chapter 2.5.7 in combination with (1) and (3) in Chapter 2.2.2]{Triebel83} we then have 
 \begin{align*}
 \holspace \b \R  \cong B^\b_{\infty,\infty}(\R) \, ,
 \end{align*}
while for open and bounded intervals $I$ we have 
\begin{align*}
B^\b_{\infty,\infty}(I) \hookrightarrow B^s_{2,2}(I)\, ,
\end{align*}
see e.g.~\cite[(7) in Chapter 3.3.1]{Triebel83}. Taking such an $I$ with $T\subset I$, the combination of both results restricted to $T$ then leads to 
\begin{align*}
 \holspace \b T  \cong B^\b_{\infty,\infty}(T) \hookrightarrow B^s_{2,2}(T) \cong H^s(T)\, ,
\end{align*}
where the last identity was noted in \eqref{eq:bessel-rkhs-onT}. Since we have already seen that almost all paths are $\b$-H\"older continuous, we thus obtain
the assertion.
%
%
\end{proofof}

\begin{proofof}{Example \ref{ex:rlp}}
 We first consider the case $\a\leq 1/2$. Here we recall from \cite[Theorem 5.4]{Picard11a} that the RKHS $H_{R^{(\a)}}$ of the Riemann-Liouville process
 and the RKHS  $H_{B^{(\a)}}$ of the fractional Brownian motion are equal up to equivalent norms. Consequently, they are essentially 
 equivalent \changebegin and we have 
 \begin{align*}
 H_{B^{(1/2)}} \subset  H_{B^{(\a)}} \smspac H_{R^{(\a)}}\, .
 \end{align*}
 Let $J:H_{B^{(1/2)}} \to H_{R^{(\a)}}$ be the corresponding  embedding map, 
 which is bounded by e.g.~\cite[Lemma 23]{ScSt25a}. 
 As in the proof of Example  \ref{ex:fbm},
 we now assume that there was an  RKHS $H$ on $T$  with  
$\P(\{\sppath R^{(\a)} \in H\}) =1$. We then choose $T^*$ as in 
the proof of Example  \ref{ex:fbm}, so that 
  Corollary \ref{cor:restrict-to-bounded} again ensures
that 
\begin{align*}
I_{(H_{R^{(\a)}})_{|T^*},\lb}:(H_{R^{(\a)}})_{|T^*}\to \Lx 2 {T^*}
\end{align*}
 is nuclear.
Since $J$ is bounded, a diagram analogous to that in
the proof of Example  \ref{ex:fbm} then shows that 
$I_{(H_{1/2})_{\changed{|T^*} },\lb}$ is nuclear. However, we have already seen 
  in the proof of 
  Example  \ref{ex:fbm} that is false.
%
 \changeend
 

 Let us now consider the case $\a\in (1/2,1)$. Here, we need to recall that we can jointly realize the processes $R^{(\a)}$ and $B^{(\a)}$
 such that the paths of $R^{(\a)} - B^{(\a)}$ are $C^\infty$, see e.g.~\cite[Theorem 5.1]{Picard11a}. 
 Since for $s\in (1/2,\a)$ all
   $C^\infty$-functions are contained in 
 $B_{2,2}^s(T)$, see e.g.~\cite[Paragraph 7.32]{AdFo03}, the identity $B_{2,2}^s(T) \cong H^s(T)$, see \eqref{eq:bessel-rkhs-onT},
 in combination with 
 \begin{align*}
  \P(\{\sppath B^{(\a)} \in H^s(T)\}) =1 \, ,
 \end{align*}
 which has been shown in Example \ref{ex:fbm}, gives $\P(\{\sppath R^{(\a)} \in H^s(T)\}) =1$.
\end{proofof}

\begin{proofof}{Example \ref{ex:matern}}
Since $T$ is bounded, there exists an open  Euclidean ball $B$ with $T\subset B$.
Moreover,
the RKHS $H_{\a,\s}(B)$ of the Mat\'ern kernel  on $B$   satisfies $H_{\a,\s}(B) \cong B_{2,2}^{\a + d/2}(B)$, see 
e.g.~\cite[Corollary 10.13]{Wendland05} in combination with e.g.~\cite[(4.15)]{RaWi06}, or alternatively 
\cite[(27)]{PoBeScOa24a}.
In addition, we have  $B_{2,2}^{\a + d/2}(B) \cong H^{\a + d/2}(B)$, see \eqref{eq:bessel-rkhs-onT}, 
and hence we find 
\begin{align*}
H_{\a,\s}(B) \cong  H^{\a + d/2}(B)\, .
\end{align*}
Restricting both spaces onto $T$ leads to $H_{\a,\s}(T) \cong  H^{\a + d/2}(T)$ by
 \cite[Proposition 3]{ScSt25a}. 
By Theorem \ref{thm:sobol-main} and $H_X = H_{\a,\s}(T)$  we then obtain the assertion.
%
%
%
%
\end{proofof}

\begin{proofof}{Example \ref{ex:abel}}
 We first note that we have $k_X(t,t) = k_X(t-t,0) = k_X(0,0)$ for all $t\in T$. Consequently, $k_X$ is bounded, see e.g.~\cite[(4.14)]{StCh08}.
 In addition, for $\kappa (r) := k_X(r,0)$ we have $k_X(t_1,t_2) = k_X(t_2,t_1) = k_X(-t_1+t_2,0) = \k(-t_1+t_2)$, and hence 
 $k_X$ is translation invariant. Finally, recall that $k_X$ is  measurable by assumption. 
 With the help of these properties of $k_X$ \cite[Lemma 4.1 and 4.2]{StZi21a} then show that $k_X$ is continuous and 
 there exist an at most countable family $(e_i^*)_{i\in I}$ of continuous functions $T\to \R$ with 
 $\sup_{i\in I}\inorm{e_i^*} \leq \sqrt 2$
 and a summable
 family $(\mu_i)_{i\in I}\subset [0,\infty)$ such that $(\eqclass {e_i^*})_{i\in I}$
 is an ONB of $\Lx 2 \lb$ and 
 \begin{align}\label{ex:abel-h1}
 k_X(t_1,t_2) = \sum_{\mu_i>0} \mu_i e_i^*(t_1) e_i^*(t_2)\, , \myqquad t_1,t_2\in T,
 \end{align}
 where the series converges absolutely for all $t_1,t_2\in T$ as well as uniformly in $(t_1,t_2)$.
 Moreover, \cite[Lemma 4.2]{StZi21a} shows that 
 the positive $\mu_i$ are the collection of  non-zero eigenvalues of $T_{k,\lb}$ including geometric multiplicities and 
 the $\eqclass {e_i^*}$ are corresponding eigenfunctions. Finally, 
for the functions defined in \eqref{eq:def-ei} with $f_i :=  {e_i^*}$ and $\mu_i>0$ we have 
 \begin{align*}
e_i = \mu_i^{-1} I_{X,\nu}^* \eqclass{e_i^*} 
= \mu_i^{-1} \sum_{\mu_j>0} \mu_j \skprodb {\eqclass{e_i^*} }{\eqclass{e_j^*} }_{\Lx 2 \lb} e_j^* = e_i^*\, ,
 \end{align*}
 where in the second   step we used \cite[Theorem 2.11]{StSc12a} in combination with \eqref{ex:abel-h1}.
 
In addition, since $k_X$ is continuous, every $h\in H_X$ is continuous. 
Moreover,  the Haar measure satisfies 
 $\lb(O)>0$ for all non-empty and open $O\subset G$. 
 Together, these properties show that 
 $\lb$ is $H_X$-positive.
 
 Using all these properties, the characterization  now follows by combining 
 \changebegin Corollary \ref{cor:general-result} \changeend
 and  Theorem  \ref{thm:ew-char}
 as discussed after Theorem  \ref{thm:ew-char}.
In addition, $\sup_{i\in I}\inorm{e_i^*} \leq \sqrt 2$ in combination with $e_i^*\in C(T)$ 
and the summability assumption on the eigenvalues 
ensures that \eqref{eq:beta-power-kx-full}
converges absolutely in $C(T\times T)$, and thus $k_X^{1/2}$ is continuous. This in turn implies $H_X^{1/2} \subset C(T)$.
\end{proofof}